\documentclass[10pt,a4paper]{article}
\usepackage{amstext}
\usepackage{array}
\usepackage{amsfonts}
\usepackage{amssymb, bm}
\usepackage{graphicx}
\usepackage{epstopdf}
\usepackage{algorithm}     %serve per scrivere lo pseudo codice
\usepackage{algpseudocode} %serve per scrivere lo pseudo codice
\usepackage{url}
\usepackage{caption}
\usepackage{subfig}
\usepackage{csquotes}
\usepackage{tikz}
\usetikzlibrary{arrows}
\usetikzlibrary{decorations.markings}
\usetikzlibrary{shapes.geometric}
\usetikzlibrary{shapes.arrows}
\usetikzlibrary{fit,calc}
\usetikzlibrary{quotes,angles}
\usetikzlibrary{positioning, arrows.meta}
\usepackage{animate}  
\usepackage{graphicx}
\usepackage{pgf}
\usepackage{bbm}
\usepackage{pgfplots}
\usepackage{mathtools}
\usepackage{physics}
\usepackage{amsthm}
\usepackage{hyperref}
\usepackage{stmaryrd}
\newcommand{\vertiii}[1]{{\left\vert\kern-0.25ex\left\vert\kern-0.25ex\left\vert #1 
   \right\vert\kern-0.25ex\right\vert\kern-0.25ex\right\vert}}

\newtheorem{theorem}{Theorem}[section]

\usepackage{geometry}
 \usepackage{fancyhdr}  
  \fancypagestyle{plain}{
\fancyhead{}
\fancyhead[C]{\hfill Submitted to Comptes Rendus Mathématique, December 2023}
 }

\title{Point-set registration in bounded domains  via the Fokker-Planck equation}
\author{Angelo Iollo, Tommaso Taddei}
\date{\today}

\begin{document}

\maketitle

\begin{abstract}
We present  a point set registration method in  bounded domains based on the solution to the Fokker Planck equation. Our approach leverages
(i) density estimation based on Gaussian mixture models;
(ii) a stabilized finite element discretization of the Fokker Planck equation;
(iii) a specialized method for the integration of the particles.
We review relevant properties of the Fokker Planck equation that provide the foundations for the numerical method.
We discuss two strategies for the integration of the particles and we propose a regularization technique to control the distance of the particles from the boundary of the domain.
We perform extensive numerical experiments for two two-dimensional model problems to illustrate the many features of the method.
\end{abstract}

\section{Introduction}
\label{sec:introduction}

Point set registration (PSR) is the process of finding a spatial transformation that aligns two point clouds. Given the domain $\Omega \subset \mathbb{R}^n$ and the reference and the point clouds 
$\{X_i^0\}_{i=1}^{N}, \{ X_j^\infty\}_{j=1}^{M} \subset \Omega$, we seek a transformation $\Phi:\Omega \to \Omega$ such that the mapped point cloud $\{ \Phi(X_i^0)  \}_{i=1}^N$ is close --- in the sense of Hausdorff --- to the target point cloud. For several applications, it is also important to identify a (possibly non-autonomous) flow of diffeomorphisms $\Phi: \Omega \times \mathbb{R}_+ \to \Omega$ to coherently deform the reference point cloud 
$\{X_i^0\}_{i=1}^{N}$
into 
$\{ X_j^\infty\}_{j=1}^{M}$. In this note,  we focus on the problem of PSR in moderate-dimensional ($n=2$ or $n=3$)  bounded domains; we target applications in model order reduction (MOR) where coordinate transformations are employed to align coherent features of the solution field \cite{cucchiara2024model,iollo2014advection,iollo2022mapping,taddei2020registration};
we also envision applications to data assimilation tasks that involve the interpolation between point clouds in complex geometries.

We propose a PSR method based on the solution to the Fokker Planck (FP) equation. First, we model the reference and the target point clouds as independent identically distributed samples from two continuous distributions with densities $\rho_0$ and $\rho_\infty$, respectively. Second, we estimate $\rho_0$ and $\rho_\infty$ using Gaussian mixture models. 
Third, we solve a suitable FP equation to determine a flow of probability densities
$t\mapsto \rho(\cdot, t)$ that is equal to 
$\rho_0$ at time $t=0$, and converges to $\rho_\infty$ as $t\to \infty$. Fourth, we exploit the 
properties of the FP equation  to obtain the desired morphing based on the 
knowledge of $\rho(\cdot, t)$ at all times.

The FP equation describes the evolution of a probability density $\rho$ under   diffusion
and drift: it reads as a linear parabolic equation with Neumann boundary conditions on $\partial \Omega$, constant diffusion and time-independent advection. The FP equation has been extensively employed in 
Information Theory 
\cite{gray2011entropy}, stochastic processes in Physics and Chemistry
\cite{van1992stochastic}, and more recently in machine learning for generative modeling 
\cite{song2019generative}. The solution to the FP equation is tightly linked to the solution to  optimal transportation (OT) problems
\cite{santambrogio2015optimal}: however, while the FP equation is a linear parabolic equation with time-independent coefficients, the evolution equation associated to OT problems is highly nonlinear and much more difficult to handle, despite the many recent contributions to the field 
\cite{peyre2019computational}.

The outline of the note is as follows.
In section \ref{sec:mathematical_background}, 
we introduce the FP equation and we discuss relevant properties that are exploited in the remainder of the note.
In section \ref{sec:methods}, we 
present the numerical method:
first, we introduce 
the  discretization of the FP equation, which is based on the finite element method for space discretization and on the Crank Nicolson method for time integration;
second, we discuss two distinct methods for the transport of the point clouds;
third, we present a regularization approach to control the distance of the particles from the boundary;
fourth, we briefly review the problem of density estimation through Gaussian mixtures.
In section \ref{sec:numerics}, we present several numerical results in two dimensions that illustrate the performance of the method.
Section \ref{sec:conclusions} concludes the paper.

\section{Mathematical background}
\label{sec:mathematical_background}

We review select properties of the Fokker Planck (FP) equation that provide the mathematical foundation for the method. We emphasize that the FP equation has been the subject of extensive studies in Analysis and Calculus of Variations; in particular, the material of this section relies heavily on 
\cite{gray2011entropy,jordan1998variational,ouhabaz2009analysis,santambrogio2015optimal}. Below, $\Omega$ is a Lipschitz domain in $\mathbb{R}^n$, while $\mathcal{P}_2(\Omega)$ refers to the space of probability density functions (pdfs) with bounded second-order moments; the acronym ``a.s.'' 
(resp., ``a.e.'')
stands for ``almost surely'' 
(resp., ``almost every'')
with respect to the Lebesgue measure in $\Omega$ and the acronym ``i.i.d.'' stands for ``independent identically distributed''.

\subsection{Kullback-Leibler   divergence}
Given the densities $\rho,\rho_\infty\in \mathcal{P}_2(\Omega)$, we introduce the Kullback-Leibler (KL) divergence between $\rho$ and $\rho_\infty$ as \cite{kullback1997information}
\begin{equation}
\label{eq:KLdivergence}
D_{\text{KL}}(\rho \parallel \rho_\infty) = \int_{\Omega} \rho(x) \log \left( \frac{\rho(x)}{\rho_\infty(x)} \right) dx.
\end{equation}
Note that $D_{\text{KL}}(\rho \parallel \rho_\infty)$ is infinite  if  the support of $\rho_\infty$ does not contain the support of $\rho$.
Theorem \ref{th:KL} summarizes key properties of the KL divergence.

\begin{theorem}
\label{th:KL}
The following hold.
\begin{enumerate}
\item 
Given the densities $\rho,\rho_\infty\in \mathcal{P}_2(\Omega)$,
$D_{\text{KL}}(\rho \parallel \rho_\infty) \geq 0$ and 
$D_{\text{KL}}(\rho \parallel \rho_\infty) =0$ if and only if $\rho=\rho_\infty$ a.s..
\item 
The KL divergence 
\eqref{eq:KLdivergence}
is asymmetric, that is
in general 
$D_{\text{KL}}(\rho \parallel \rho_\infty) \neq 
D_{\text{KL}}(\rho_\infty \parallel \rho)$.
\item 
The KL divergence 
\eqref{eq:KLdivergence}
is strongly convex, that  is, 
given $\rho_0,\rho_1,\rho_\infty\in \mathcal{P}_2(\Omega)$ and $\lambda\in (0,1)$, we have
$$
D_{\text{KL}}( \lambda \rho_0 + (1-\lambda) \rho_1 \parallel \rho_\infty)
\leq
 \lambda
D_{\text{KL}}( \rho_0  \parallel \rho_\infty)
+
(1-\lambda)
D_{\text{KL}}(\rho_1 \parallel \rho_\infty),
$$
where the latter holds with equality if and only if $\rho_0=\rho_1$ a.s..
\item 
(Pinsker's inequality)
The KL divergence 
\eqref{eq:KLdivergence} satisfies
$
\|\rho  - \rho_\infty  \|_{L^1(\Omega)}
\leq
\sqrt{
2 D_{\text{KL}}( \rho  \parallel \rho_\infty)
}$.
\end{enumerate}
\end{theorem}

We provide the proof of Theorem \ref{th:KL} in Appendix \ref{sec:proofs}.  We notice that the KL divergence is minimized when $\rho=\rho_\infty$ (cf. Property 1); as discussed below, the strong convexity property and the Pinsker's inequality
are the key to prove the existence of paths that link the initial and the target distributions (cf. Theorem \ref{th:FPequation}). We further remark that the KL divergence is not a distance in $\mathcal{P}_2(\Omega)$ since it is not symmetric (cf. Property 2). If we denote by $V=-\log(\rho_\infty)$ the potential associated with the density $\rho_\infty$, we can express \eqref{eq:KLdivergence} as:
\begin{equation}
\label{eq:KLdivergence_revisited}
D_{\text{KL}}(\rho \parallel \rho_\infty) = \int_{\Omega} F[\rho] \, dx
=
\underbrace{\int_\Omega \rho \log (\rho) \, dx}_{\rm =: (I)}
\, + \,
\underbrace{\int_\Omega \rho V \, dx}_{\rm =: (II)},
\quad
{\rm where} \quad F[\rho] = 
 \rho \log (\rho) + \rho V.
\end{equation}
In Bayesian inference,  the first term in \eqref{eq:KLdivergence_revisited} is interpreted as an entropy function, while the second term is interpreted as a potential energy and the sum of the two terms is interpreted as a free-energy function (see, e.g., 
\cite{fox2012tutorial}).

\subsection{Fokker-Plank equation}

Let $\rho, \rho_\infty \in \mathcal{P}_2(\Omega)$ and let 
$V=-\log(\rho_\infty)$; then, we define the Fokker-Plank (FP) equation:
\begin{equation}
\label{eq:FPequation}
\left\{
\begin{array}{ll}
\displaystyle{
\partial_t \rho + \nabla \cdot \left( - \rho \nabla V - \nabla \rho 
\right) = 0
}     &  \displaystyle{{\rm in} \; \Omega \times \mathbb{R}_+}\\[3mm]
\rho(\cdot,0) = \rho_0     & 
{\rm in} \; \Omega  \\[3mm]
\left( \rho \nabla V + \nabla \rho 
\right) \cdot \mathbf{n} = 0     & 
{\rm on} \; \partial \Omega \times \mathbb{R}_+  \\
\end{array}
\right.
\end{equation}
where $\mathbf{n}$ denotes the outward normal to $\partial \Omega$.
We assume that $\rho_\infty$ is strictly positive in $\Omega$ and is sufficiently smooth to ensure that $\nabla V$ is Lipschitz continuous. 
We introduce the spaces $\mathcal{H}=L^2(\Omega)$ and
$\mathcal{V}=H^1(\Omega)$ and the   weak formulation of the FP equation \eqref{eq:FPequation}:
find
$\rho \in C(\mathbb{R}_+; \mathcal{H})$ such that $\rho(0) = \rho_0$ and 
\begin{equation}
\label{eq:weak_form}
\int_\Omega  \left( \partial_t \rho(t) \cdot v
+ ( \nabla \rho(t) + \rho(t) \nabla V ) \cdot \nabla v  \right) \, dx = 0
\quad
\forall \, v\in \mathcal{V}, 
\quad
{\rm a.e.} \;\; t\in \mathbb{R}_+.
\end{equation}

Theorem \ref{th:FPequation} reviews relevant properties of the FP equation; the proof is postponed to Appendix \ref{sec:proofs}.

\begin{theorem}
\label{th:FPequation}
Let $\Omega$ be a Lipschitz domain in $\mathbb{R}^n$; let 
$\rho_0, \rho_\infty \in \mathcal{P}_2(\Omega)$ satisfy
(i)
$\inf_{x\in \Omega} \rho_\infty (x) > 0$,
and 
(ii) $\nabla V $  is Lipschitz-continuous in $\Omega$.
The following hold.
\begin{enumerate}
\item
There exists a unique weak solution $\rho \in C(\mathbb{R}_+; \mathcal{H})$ to \eqref{eq:FPequation}.
\item
The solution  to \eqref{eq:FPequation} satisfies
$\rho(\cdot,t) \in \mathcal{P}_2(\Omega)$ for all $t>0$.
\item
The solution to \eqref{eq:FPequation} converges to $\rho_\infty$ for $t\to \infty$ in $L^1$, that is,
$\lim_{t\to \infty} \|\rho(\cdot, t) - \rho_\infty  \|_{L^1(\Omega)} = 0$.
\item
Let $u = - \nabla (\log (\rho) + V)$ and let $X_1^0,\ldots,X_N^0 \overset{\rm iid}{\sim} \rho_0$. Consider the trajectories
\begin{equation}
\label{eq:trajectories_FP}
\left\{
\begin{array}{ll}
\dot{X}_i(t) = u(X_i(t), t) & t>0, \\
X_i(0) = X_i^0, & \\
\end{array}
\right.
\quad
i=1,\ldots,N.
\end{equation}
Then, 
for any $t>0$, 
we have $X_1(t),\ldots,X_N(t) \overset{\rm iid}{\sim} \rho(\cdot, t)$.
\end{enumerate}
\end{theorem}
Theorem \ref{th:FPequation} provides the foundations for the use of the FP equation for point-set registration. The first statement ensures the existence and the uniqueness of the solution to \eqref{eq:FPequation};  the second statement ensures that $\rho(\cdot, t)$ is a pdf for all $t>0$; the third statement shows that $\rho(\cdot, t)$ converges to the target distribution for $t\to \infty$; finally, the fourth statement provides a constructive way to sample from the target distribution and also to generate trajectories that link initial and target pdfs.

\subsection{Fokker-Plank equation and gradient flows}
\label{sec:gradient_flows}
The authors of \cite{jordan1998variational} showed that the FP equation can be interpreted as the gradient flow of the KL divergence with respect to the Wasserstein metric.
We can hence interpret the solution to the time-discrete FP equation as the solution to a sequence of (linearized) optimal transport problems. Below, we provide an informal proof of this statement for $\Omega=\mathbb{R}^n$; we refer to \cite{jordan1998variational} for the rigorous derivation and for further details.  

First, we recall the 2-Wasserstein metric between the densities $\rho_0$ and $\rho_1$:
\begin{equation}
\label{eq:wasserstein_metric}
W_2^2(\rho_0,\rho_1) :=\inf_{X:\mathbb{R}^n\to \mathbb{R}^n} \int_{\mathbb{R}^n} \big| X(\xi) - \xi \big|^2 \rho_0(\xi) \, d\xi\quad
{\rm s.t} \;\; \rho_1(X(\xi)) {\rm det} \nabla X(\xi) = \rho_0(\xi).    
\end{equation}
It is possible to show that the optimal transport map $X$ is of the form $X = \texttt{id} + \nabla \psi $ where $\texttt{id}(\xi) = \xi$ is the identity map and $ \psi : \mathbb{R}^n \to \mathbb{R}$ (cf. \cite{brenier1991polar}, see also \cite[Chapter 1.3.1]{santambrogio2015optimal}). Therefore, we can restate \eqref{eq:wasserstein_metric} as
\begin{equation}
\label{eq:wasserstein_metric_step2}
W_2^2(\rho_0,\rho_1) :=\inf_{\psi:\mathbb{R}^n\to \mathbb{R}} \int_{\mathbb{R}^n} \big| \nabla \psi(\xi) \big|^2 \rho_0(\xi) \, d\xi\quad
{\rm s.t} \;\; 
\rho_1(\texttt{id} + \nabla \psi(\xi)) {\rm det} \left(  \nabla^2 \psi(\xi) \right) = \rho_0(\xi).    
\end{equation}

%\rho_1(X(\xi)) {\rm det} \nabla X(\xi) = \rho_0(\xi).    

We denote by $\{t^{(k)} = k \Delta t\}_{k \in \mathbb{N}}$ the time grid, and we denote by $\rho^{(k)}$ the estimate of $\rho$ at time $t^{(k)}$; we assume that $\rho^{(k+1)}$ results from the evaluation of $\rho^{(k)}$ subject to a potential velocity field $\epsilon \nabla \psi$ with $\epsilon \ll 1$. We express
$\rho^{(k+1)} = \rho^{(k)} + \epsilon \eta$. Since 
$\rho^{(k+1)}$ satisfies
$\rho^{(k+1)}(X(\xi)) {\rm det} \nabla X(\xi) = \rho^{(k)}(\xi)$ with $X = \texttt{id} + \epsilon \nabla \psi$, exploiting the Jacobi formula for the derivative of the determinant, we find
\begin{equation}
\label{eq:expression_eta}
\eta = - \nabla \cdot (\rho^{(k)} \nabla \psi) + \mathcal{O}(\epsilon).
\end{equation}

We seek $\psi$ to maximize the first variation of the KL divergence for a fixed Wasserstein distance:
\begin{equation}
\label{eq:wasserstein_gradient}
\psi^{\star} : = {\rm arg} \inf_{\psi} \int_{\mathbb{R}^n} F'[\rho^{(k)}] \cdot \eta(\psi) \, d x
\quad
{\rm s.t.} 
\left\{
\begin{array}{l}
\displaystyle{
\int_{\mathbb{R}^n} \big| \nabla \psi  \big|^2 \rho^{(k)}\, dx
=c,
} \\
\eta(\psi) =     - \nabla \cdot (\rho^{(k)} \nabla \psi) ,\\  
\end{array}
\right.
\end{equation}
where $c$ is a given constant and 
$F'[\rho ] = \log(\rho) + 1 + V $.
If we integrate by part the objective of \eqref{eq:wasserstein_gradient}, we find
$$
\int_{\mathbb{R}^n}
F'[\rho^{(k)}] \cdot \eta(\psi) \, d x
=
-
\int_{\mathbb{R}^n}
\nabla \cdot \left(
\rho^{(k)}
F'[\rho^{(k)}]
\right) \psi
 \, d x
$$
In conclusion, we can restate \eqref{eq:wasserstein_gradient} as
\begin{equation}
\label{eq:wasserstein_metric_step3}
\psi^{\star} : = {\rm arg} \inf_{\psi} -
\int_{\mathbb{R}^n}
\nabla \cdot \left(
\rho^{(k)}
\nabla
F'[\rho^{(k)}]
\right) \psi
 \, d \xi
\quad
{\rm s.t.} 
\quad
\int_{\mathbb{R}^n} \big| \nabla \psi  \big|^2 \rho^{(k)} \, d x
=c.
\end{equation}

The optimality conditions for \eqref{eq:wasserstein_metric_step3}  are
\begin{equation}
\label{eq:wasserstein_metric_step4}
\left\{
\begin{array}{ll}
\displaystyle{
- \int_{\mathbb{R}^n} \nabla \cdot 
\left(    \rho^{(k)}
\nabla
F'[\rho^{(k)}]     \right) \cdot \delta \psi \, dx
+  \lambda \int_{\mathbb{R}^n} 
\rho^{(k)}  \nabla  \psi  \cdot 
 \nabla   \delta \psi \, dx
 = 0
}
&
\forall \, \delta \psi \in  \mathcal{V} ;
\\[3mm]
\displaystyle{
\int_{{\mathbb{R}^n}} \big| \nabla \psi \big|^2 \rho^{(k)} \, d x 
=c;
}
&
\\
\end{array}
\right.
\end{equation}
where $\mathcal{V}:=  H^1({\mathbb{R}^n})$.
If we integrate by part the second term in the first equation of \eqref{eq:wasserstein_metric_step4}, we find\footnote{We here exploit the fact that the space $\mathcal{V}$ is dense in $L^2$.}
$$
\nabla \cdot 
\left(
\rho^{(k)} \nabla \psi
\right)
=
-\frac{1}{\lambda} 
\nabla \cdot 
\left(
\rho^{(k)}  F'[ \rho^{(k)}   ]
\right)
$$
where the constant $\lambda$ is a function of $c$ that is determined by substituting the latter expression in \eqref{eq:wasserstein_metric_step4}$_2$.
We conclude that for an appropriate choice of $c$ in \eqref{eq:wasserstein_gradient} we have
$$
\rho^{(k+1)}
=
\rho^{(k)}
-
\Delta t
\nabla \cdot \left(
\rho^{(k)} \nabla V + 
\nabla  \rho^{(k)},
\right)
$$
which is the explicit Euler discretization of the FP equation.

\subsection{Transport of Gaussian random variables}
We denote by $\mathcal{N}(\cdot, \mu, \sigma^2)$ the pdf associated with an univariate Gaussian random variable with mean $\mu$ and variance $\sigma^2$. We assume that 
$\rho_0:=\mathcal{N}(\cdot, \mu, \gamma^2)$ and
$\rho_\infty=\mathcal{N}(\cdot, 0, \sigma^2)$.
It is possible to show that the solution to the FP equation in $\Omega=\mathbb{R}$ is given by
(see, e.g., \cite[Chapter VIII.4]{van1992stochastic})
\begin{equation}
\label{eq:FP_gaussian}
\rho^{\rm fp}(x,t)
=
\mathcal{N}(x , \mu(t), \Sigma(t)),
\quad
{\rm where} \;\;
\mu(t) = e^{-
\frac{t}{\sigma^2}} \, \mu,
\quad
\Sigma(t) 
= 
e^{-\frac{2 t}{\sigma^2}} 
\left(
\gamma^2
+
\sigma^2
\left(
e^{\frac{2 t}{\sigma^2}} - 1
\right)
\right).
\end{equation}
For comparison, we recall the McCann interpolation between the same two Gaussian random variables, which corresponds to the geodesic (shortest path) in the Wasserstein metric
(see, e.g., \cite[Remark 2.31]{peyre2019computational})
\begin{equation}
\label{eq:OT_gaussian}
\rho^{\rm ot}(x,t)
=
\mathcal{N}(x , \mu(t), \Sigma(t)),
\quad
{\rm where} \;\;
\mu(t) = (1- t)\mu,
\quad
\Sigma(t) 
= 
\left(
(1-t)
\gamma 
+
t
\sigma 
\right)^2,
\quad
t\in [0,1].
\end{equation}

We observe that the solution to \eqref{eq:FPequation} connects $\rho_0$ to $\rho_\infty$ in infinite time at non-constant speed, while the optimal transport map achieves  the same goal in finite time and at constant speed. Note, nevertheless, that the convergence of $\rho^{\rm fp}$ to $\rho_\infty$ is exponentially-fast: this observation motivates the use of adaptive time stepping strategies to integrate \eqref{eq:FEM_weak}.

\section{Methodology}
\label{sec:methods}

\subsection{High-fidelity discretization}
\label{sec:FEM}
We rely on the finite element (FE) method to discretize \eqref{eq:weak_form} in space, and on the Crank-Nicolson method for time integration. Given the time grid $\{t^{(k)} = k \Delta t  \}_{k\in \mathbb{N}}$ and the FE space $\widehat{\mathcal{V}} \subset \mathcal{V}$, we define the sequence of approximate solutions $\{ \widehat{\rho}^{(k)}  \}_{k\in \mathbb{N}} \subset \widehat{\mathcal{V}}$ such that
$\widehat{\rho}^{(0)}= \rho_0$ and
\begin{equation}
\label{eq:FEM_weak}
\int_\Omega  \left( 
\frac{\widehat{\rho}^{(k+1)} - \widehat{\rho}^{(k)}}{\Delta t}
\right) 
  \cdot v
+ \frac{1}{2} \left( 
\nabla \widehat{\rho}^{(k+1)}+  \widehat{\rho}^{(k+1)} \nabla V 
+
\nabla \widehat{\rho}^{(k)}+  \widehat{\rho}^{(k)} \nabla V 
\right) 
\cdot \nabla v    \, dx = 0
\quad
\forall \, v\in \widehat{\mathcal{V}}, 
\end{equation}
for $k=0,1,\ldots$.
We  remark  that 
problem \eqref{eq:FEM_weak} might feature strong advection: for this reason, 
in the numerical experiments, 
we resort to  a SUPG stabilization of \eqref{eq:FEM_weak} to avoid spurious oscillations \cite{brooks1982streamline}.

Given the particles
$\{ X_j^{(0)} \}_{j=1}^N$, we consider two distinct strategies
to generate the approximate trajectories:
the \emph{``ODE'' method} and the 
\emph{``gradient flow (GF)'' method}.
The former is justified by Property 4 in in Theorem \ref{th:FPequation}, while the 
latter is motivated by the discussion in section \ref{sec:gradient_flows}.

\paragraph{\textbf{ODE method.}}
In this approach, we 
discretize 
the ODE system
\begin{equation}
\label{eq:ODE_refined}
\left\{
\begin{array}{ll}
\dot{X}_i(\tau) = \widehat{\mathbf{u}}({X}_i(\tau), \tau),     & 
  \tau \in (t^{(k)}, t^{(k+1)}) ,\\[3mm]
{X}_i(t^{(k)})
= {X}_i^{(k)},& \\
\end{array}
\right.    
\end{equation}
where $\widehat{\mathbf{u}}(\cdot, \tau)   = \nabla V + \nabla \widehat{\rho}(\cdot , \tau)  $ and $\widehat{\rho}(\cdot , \tau)$ is a suitable approximation of the density in the interval $(t^{(k)}, t^{(k+1)})$. In the numerical experiments, we 
consider an explicit Euler discretization of \eqref{eq:ODE_refined} based on a single time step:
\begin{equation}
\label{eq:explicit_euler}
X_i^{(k+1)}
=
X_i^{(k)}
+
\Delta t \, 
\widehat{\mathbf{u}}^{(k)}(X_i^{(k)}),
\quad
{\rm where} \;\; 
\widehat{\mathbf{u}}^{(k)} := \nabla V + \nabla \widehat{\rho}^{(k)},
\quad
i=1,\ldots,N, \;\;
k=0,1,\ldots.
\end{equation}
We also consider   an explicit second-order  Runge-Kutta (RK2) scheme for \eqref{eq:ODE_refined} based on the  linear interpolation of 
the density field
$$
\widehat{\rho}(\cdot , \tau)=
\left(
1 - \frac{\tau -t^{(k)} }{\Delta t}
\right)
\widehat{\rho}^{(k)} 
\, + \,
\left(
 \frac{\tau -t^{(k)} }{\Delta t}
\right)
\widehat{\rho}^{(k+1)},\quad
\tau \in 
(t^{(k)}, t^{(k+1)}).
$$
To further increase accuracy, we also lower the time step used for \eqref{eq:ODE_refined}.

\paragraph{\textbf{GF method.}}
Given the estimates $ \widehat{\rho}^{(k)}, \widehat{\rho}^{(k+1)}$ obtained through \eqref{eq:FEM_weak}, we compute the potential $\widehat{\psi}^{(k)} \in \widehat{\mathcal{V}}$ such that
\begin{subequations}
\label{eq:GF_method}
\begin{equation}
\label{eq:GF_method_a}
\int_\Omega 
\widehat{\rho}^{(k)}
\nabla \widehat{\psi}^{(k)} \cdot \nabla v \, dx
=
\int_\Omega 
\left(
\widehat{\rho}^{(k+1)} - 
\widehat{\rho}^{(k)}  \right)
 v \, dx,
 \quad
 \forall \, v\in \widehat{\mathcal{V}},
\end{equation}    
and then we estimate the positions of the particles as
\begin{equation}
X_i^{(k+1)}
=
X_i^{(k)}
+
 \, 
\widehat{\mathbf{u}}^{(k)}(X_i^{(k)}),
\quad
{\rm where} \;\; 
\widehat{\mathbf{u}}^{(k)} := \nabla \psi^{(k)},
\quad
i=1,\ldots,N, \;\;
k=0,1,\ldots.    
\end{equation}
\end{subequations}
 Note that \eqref{eq:GF_method_a} corresponds to the linearization of the Monge-Ampere equation; in future works, we  shall exploit this analogy to determine higher-order approximations of the trajectories of the particles. We observe that the linear system associated with \eqref{eq:GF_method_a}
 is poorly conditioned when $\inf_{x\in \Omega} \widehat{\rho}^{(k)}$ is small. In the numerical examples, we consider the regularized problem:
\begin{equation}
\label{eq:GF_method_regularized}
\int_\Omega 
\widehat{\rho}^{(k)}
\nabla \widehat{\psi}^{(k)} \cdot \nabla v  \,
+  \, \epsilon  \, 
\widehat{\psi}^{(k)}   \,  v
\, dx
=
\int_\Omega 
\left(
\widehat{\rho}^{(k+1)} - 
\widehat{\rho}^{(k)}  \right)
 v \, dx,
 \quad
 \forall \, v\in \widehat{\mathcal{V}},
\end{equation}    
where $\epsilon>0$ is a regularization parameter that is set equal to $10^{-10}$ in the numerical experiments.

\subsection{Treatment of boundaries}
\label{sec:bnd_regularization}
For several applications, we envision that it is necessary to control the distance of the particles from the boundaries of the domain. In our setting, this can be achieved by introducing a regularization in the potential $V$. In more detail, we consider
the regularized potential
\begin{equation}
\label{eq:regularized_potential}
V_\epsilon(x) = - \log \left( \rho_\infty(x)  \right)
\; + \; 
\frac{\epsilon}{w_\delta(x)},
\end{equation}
where $w_\delta: \Omega \to \mathbb{R}_+$ is a (regularized) distance function from the boundary of $\Omega$.
The second term in \eqref{eq:regularized_potential} 
introduces a repulsive force that is active in the proximity of the boundary and prevents particles from approaching   the boundary too closely.
The constant $\epsilon>0$ is a regularization parameter whose choice is investigated in the numerical examples, while the constant $\delta>0$ is the smoothing parameter that is introduced below.

In order to define the function $w_\delta$ in \eqref{eq:regularized_potential}, we first introduce the distance function
$w(x) = \max\big( {\rm dist}(x,\partial \Omega)$, 
$\texttt{tol} \big)$, where 
the user-defined tolerance $\texttt{tol}$  prevents division by zero in \eqref{eq:regularized_potential}.
Since the distance function $w$
is in general only Lipschitz continuous, we apply a smoothing to ensure that $\nabla V$ is sufficiently smooth.
We consider a non-conforming approximation of the biharmonic equation for smoothing \cite{mozolevski2007hp}.
We denote by 
$\texttt{F}_{j}$ the j-th facet of the mesh, by 
$\mathbf{n}^+$ the positive normal to each facet, by  
$ \texttt{D}_{k}$ the $k$-th element of the mesh, by 
$\texttt{I}_{\rm int}$ the set of internal facets of the mesh.
Given the field $w:\Omega \to \mathbb{R}$, we define the positive and the negative limits on the facets of the mesh
$$
w^\pm(x):=\lim_{\epsilon\to 0^+} w(x\mp \epsilon \mathbf{n});
$$
then, 
we introduce the 
normal jump of the gradient at each facet
$\llbracket   \nabla w  \rrbracket = \mathbf{n}^+ \cdot (\nabla w^+ -
\nabla w^-)$ and the facet average
$\{w\} = \frac{1}{2} (w^+ + w^-)$.
Exploiting the previous definitions, we define 
the  smoothing problem: find $w_\delta \in \widehat{\mathcal{V}}$ such that
$w_\delta\big|_{\partial \Omega}=\texttt{tol}$ and 
\begin{subequations}
\label{eq:smoothing_problem}
\begin{equation}
a_\delta(w_\delta, v) = 
\frac{1}{\delta} \int_\Omega w \cdot v \, dx   \quad
\forall \, v\in \widehat{\mathcal{V}}  \; : \;
v\big|_{\partial \Omega}=0;
\end{equation}
where
\begin{equation}
a_\delta(w,v) =   
\sum_{k=1}^{N_{\rm e}}
 \int_{ \texttt{D}_{k}   } \left( \Delta w  \cdot  \Delta v
 + \frac{1}{\delta} w \cdot v
 \right) \, dx
 +  \sum_{j \in \texttt{I}_{\rm int}}  \int_{ \texttt{F}_{j}   }
\left(
\beta_j 
\llbracket   \nabla w  \rrbracket \cdot 
\llbracket   \nabla v  \rrbracket
+ 
\frac{1}{\beta_j}
\{  \Delta w  \} \cdot \{  \Delta v  \}  \right) \, ds,
\end{equation}
and $\beta_j = \sigma_{\beta} \kappa^2 | \texttt{F}_{j}   |^{-1}$ with 
$ \sigma_{\beta} =10$, and $\kappa$ equal to the polynomial degree of the FE basis.
\end{subequations}

\subsection{Estimate of the target density}
\label{sec:target_density}
The application of the FP equation to PSR requires the estimate of the densities $\rho_0$ and $\rho_\infty$. Given the initial and the target point clouds $\{  X_i^{0} \}_{i=1}^N$ and
$\{  X_j^{\infty} \}_{j=1}^M$, we fit two Gaussian mixture models (GMMs) to determine the estimates of 
$\rho_0$ and $\rho_\infty$ that enter in \eqref{eq:FEM_weak} (see, e.g., 
 \cite[Chapter 8]{hastie2009elements}).
We resort to the Matlab function \texttt{fitgmdist} to determine the GMMs associated with the point clouds. We determine the number of mixture models using the Aikake information criterion (AIC); furthermore, we introduce a regularization of the covariance matrix to robustify the learning procedure\footnote{We consider the regularization value  $10^{-2}$; see the  Matlab documentation for further details.}.

The use of GMMs for density estimation ensures that $\rho_0$  and $\rho_\infty$ are strictly positive for all $x\in \Omega$, which is crucial to apply the results of Theorem \ref{th:FPequation}. Furthermore, we recall that for Gaussian distributions $\rho = \mathcal{N}(\cdot; \mu , \Sigma)$ we have
$$
\log \left(
\mathcal{N}(x; \mu , \Sigma) 
\right)
= -\frac{1}{2} (x-\mu)^\top \Sigma^{-1} (x-\mu).
$$
This implies that for GMMs the gradient of $\log(\rho)$ grows linearly with respect to the distance from the centers of the mixture models.

\section{Numerical results}
\label{sec:numerics}

\subsection{Transport of Gaussian distributions across a cylinder}
\label{sec:numerics_test1}

\subsubsection{Problem setting}
\begin{figure}[h!]
\centering
\subfloat[]{ 
\includegraphics[width=.47\textwidth]{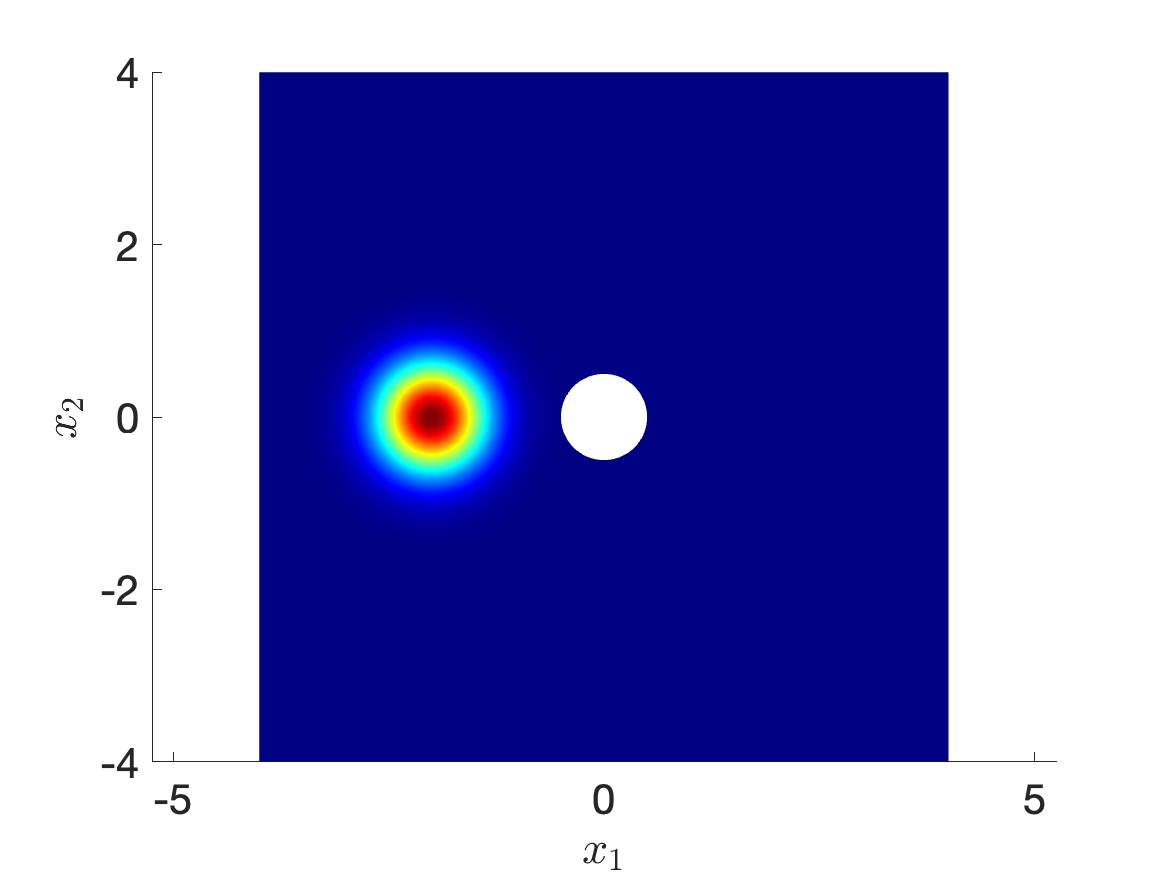}
}
~~
\subfloat[]{
\includegraphics[width=.47\textwidth]{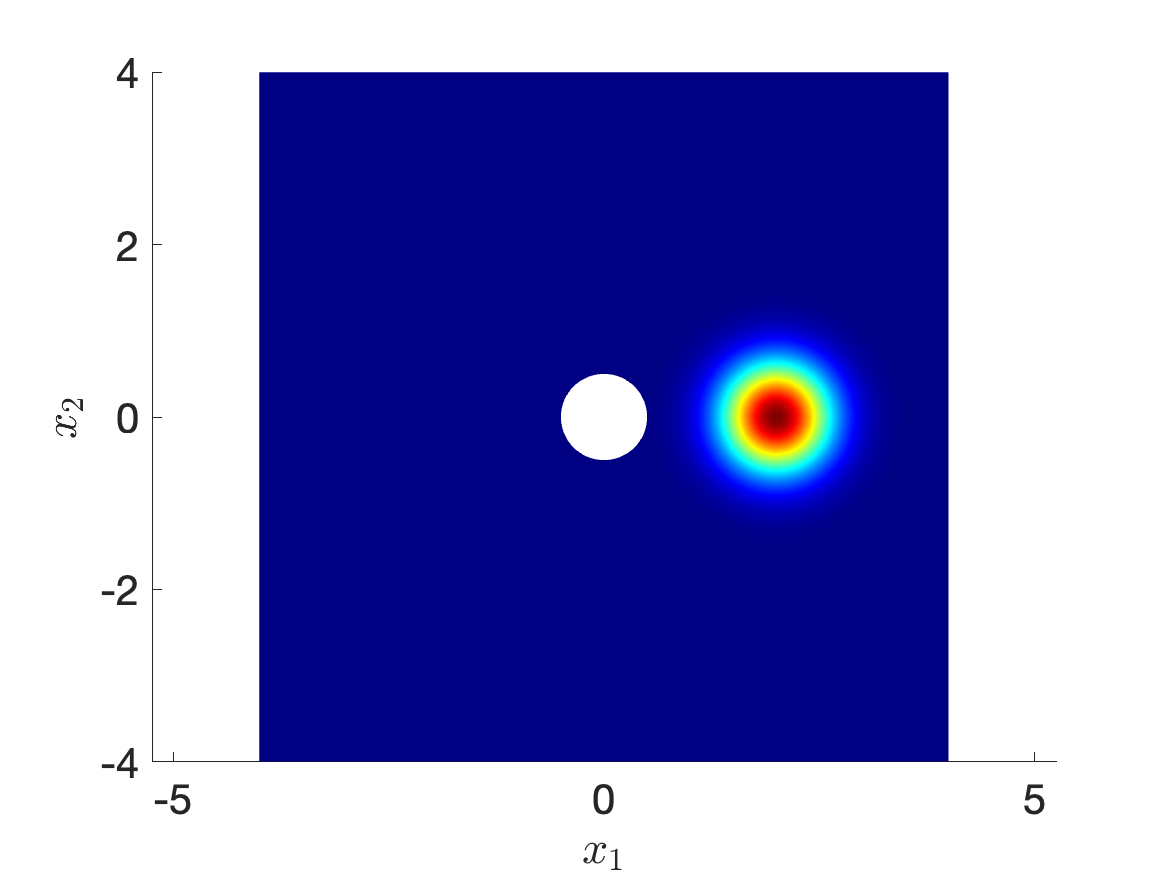}
}

\caption{transport of Gaussian distributions across a cylinder. 
(a) initial distribution.
(b) target distribution.
}
\label{fig:test1_vis}
\end{figure}

We consider the problem of transporting the Gaussian density $\rho_0 = \mathcal{N}(\cdot; \mu_0, \Sigma_0)$ with $\mu_0 = [-2,0]$ and $\Sigma_0 = 0.2 \mathbbm{1}$ to 
  $\rho_\infty = \mathcal{N}(\cdot; \mu_\infty, \Sigma_\infty)$ with $\mu_\infty = [2,0]$ and $\Sigma_\infty =  \Sigma_0$ in the domain $\Omega = (-4,4)^2 \setminus \mathcal{B}_{r=0.5}(0)$. Figure \ref{fig:test1_vis} shows the initial and the target distributions.

We resort to a P2 FE discretization with $N_{\rm hf}=6469$ degrees of freedom; we integrate  the equation in the time interval $(0,T=5)$ and we consider the grid
\begin{equation}
\label{eq:non_uniform_timegrid}
\left\{ 
t^{(k)} = T  \left( \frac{k}{K} \right)^{1.5}
\right\}_{k=0}^K, \quad
K=3000.
\end{equation}
The non-homogeneous time discretization is motivated by the observation that the dynamics is much faster at the onset of the simulation; in the future, we plan to implement more rigorous adaptive time stepping strategies. We further remark that $\nabla V$ is computed  by differentiating the FE approximation of $V= - \log(\rho_\infty)$: we notice that approximating the density $\rho_\infty$ through its nodal representation and then computing 
  $\nabla V$ as $\frac{\nabla \rho_\infty}{\rho_\infty}$ leads to significantly less accurate results. The regularized distance function in \eqref{eq:regularized_potential} is computed using $\delta=10^{-2}$ and $\texttt{tol}=10^{-2}$ and several values of $\epsilon$. 

  To validate the PSR procedure, we generate 
  $X_1^0,\ldots,X_N^0 
  \overset{\rm idd}{\sim} \rho_0$
with $N=100$ and we compare the trajectories obtained using our method with $M=100$ iid samples from 
the target distribution. Time integration of the particles is performed using the ODE method based on the explicit Euler method \eqref{eq:explicit_euler}
 with the same time step  that is used for the integration of \eqref{eq:FEM_weak}.
We also considered the GF method and the ODE method with RK2 time integration: since the results of the latter two methods are nearly equivalent to the ones obtained using 
\eqref{eq:explicit_euler}, we do not report them here.
 
\subsubsection{Results}
\begin{figure}[h!]
\centering
\subfloat[]{ 
\includegraphics[width=.45\textwidth]{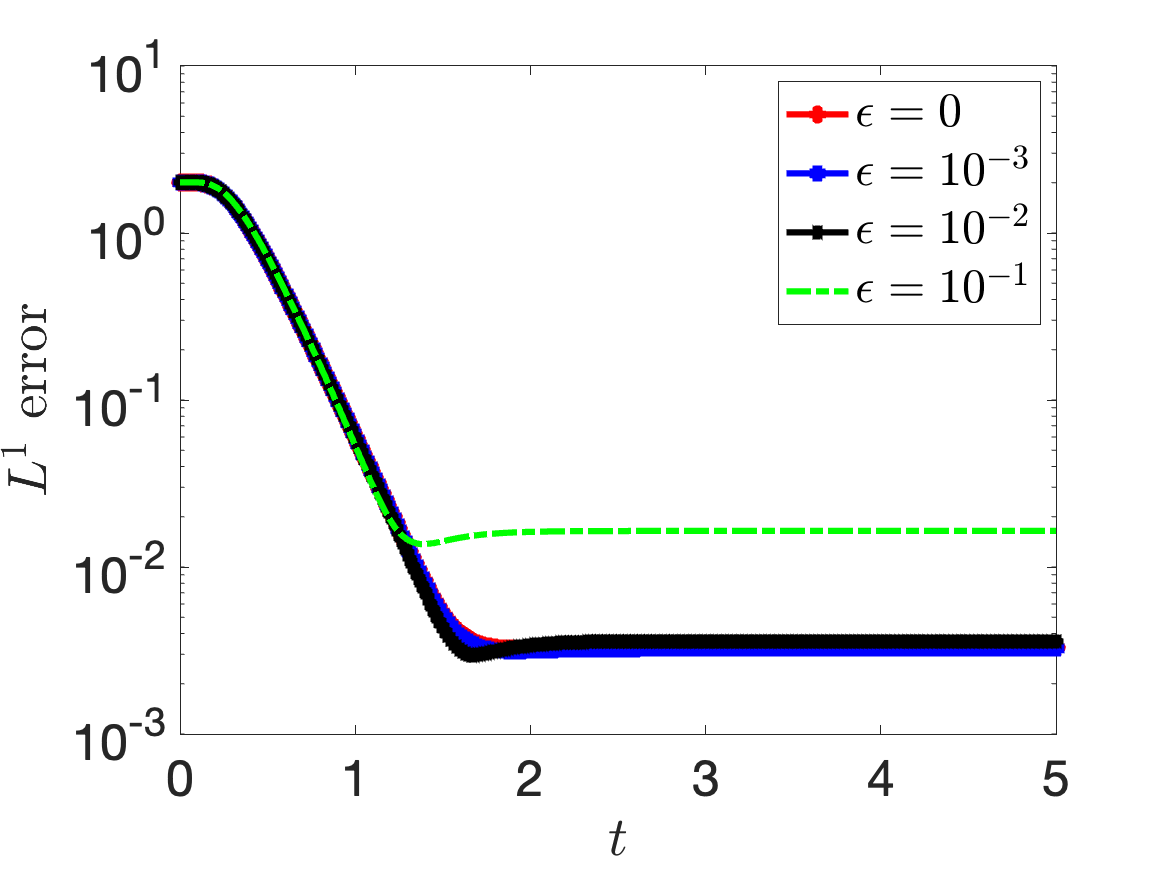}
}
~~
\subfloat[]{
\includegraphics[width=.45\textwidth]{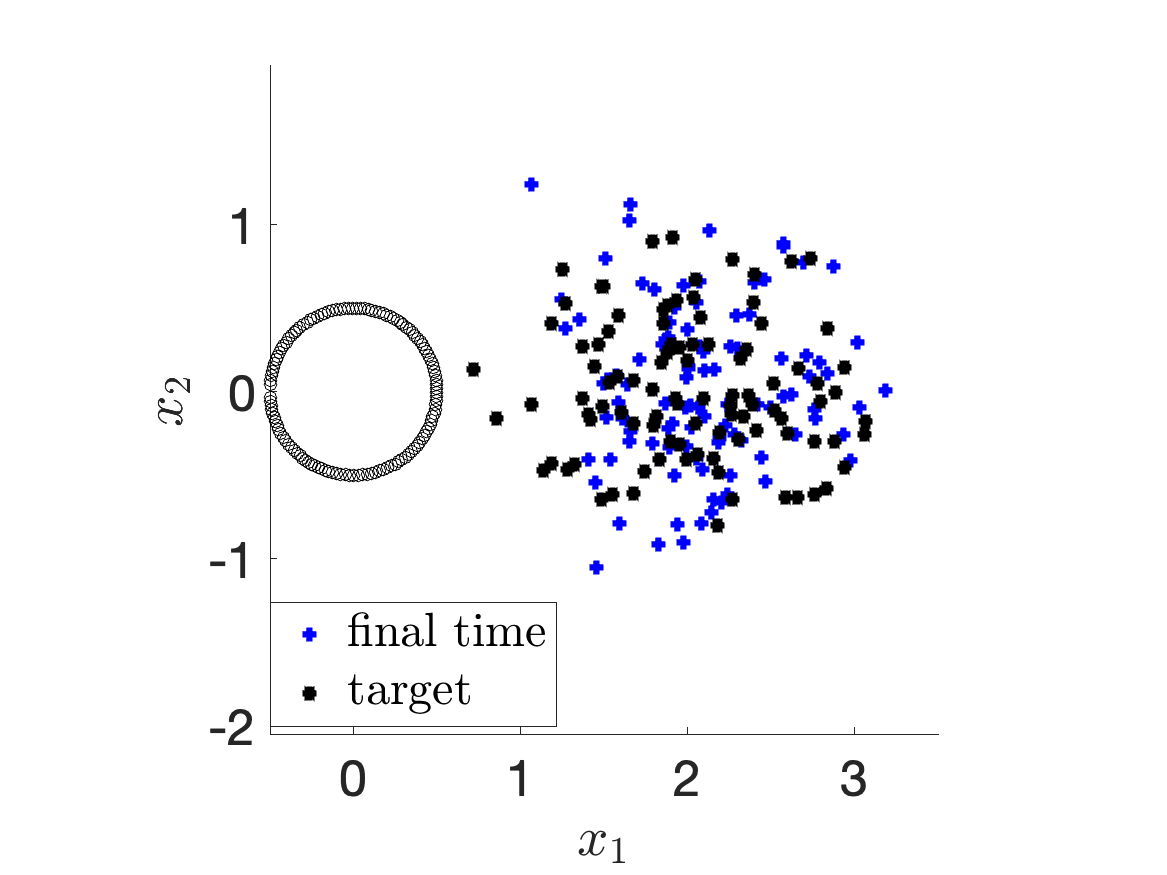}
}

\subfloat[]{
\includegraphics[width=0.45\textwidth]{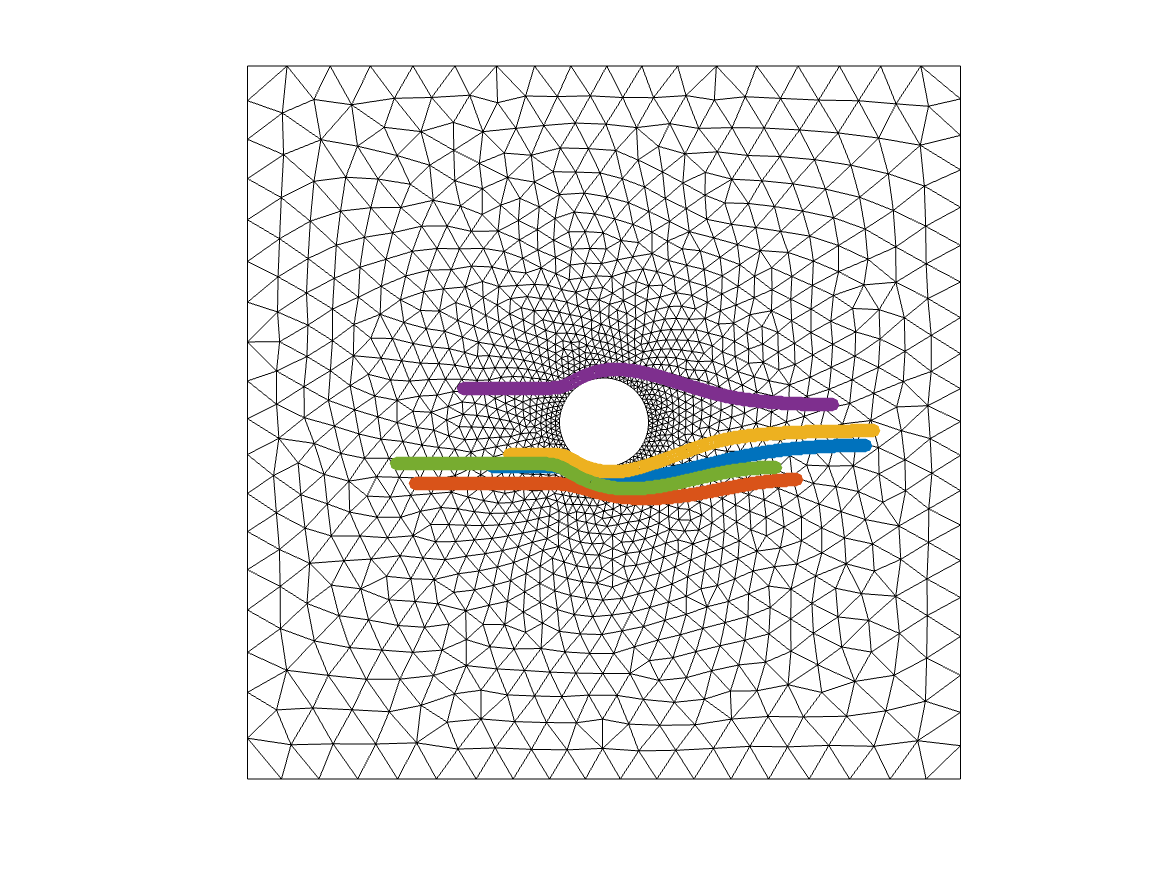}
}
~~
\subfloat[]{
\includegraphics[width=.45\textwidth]{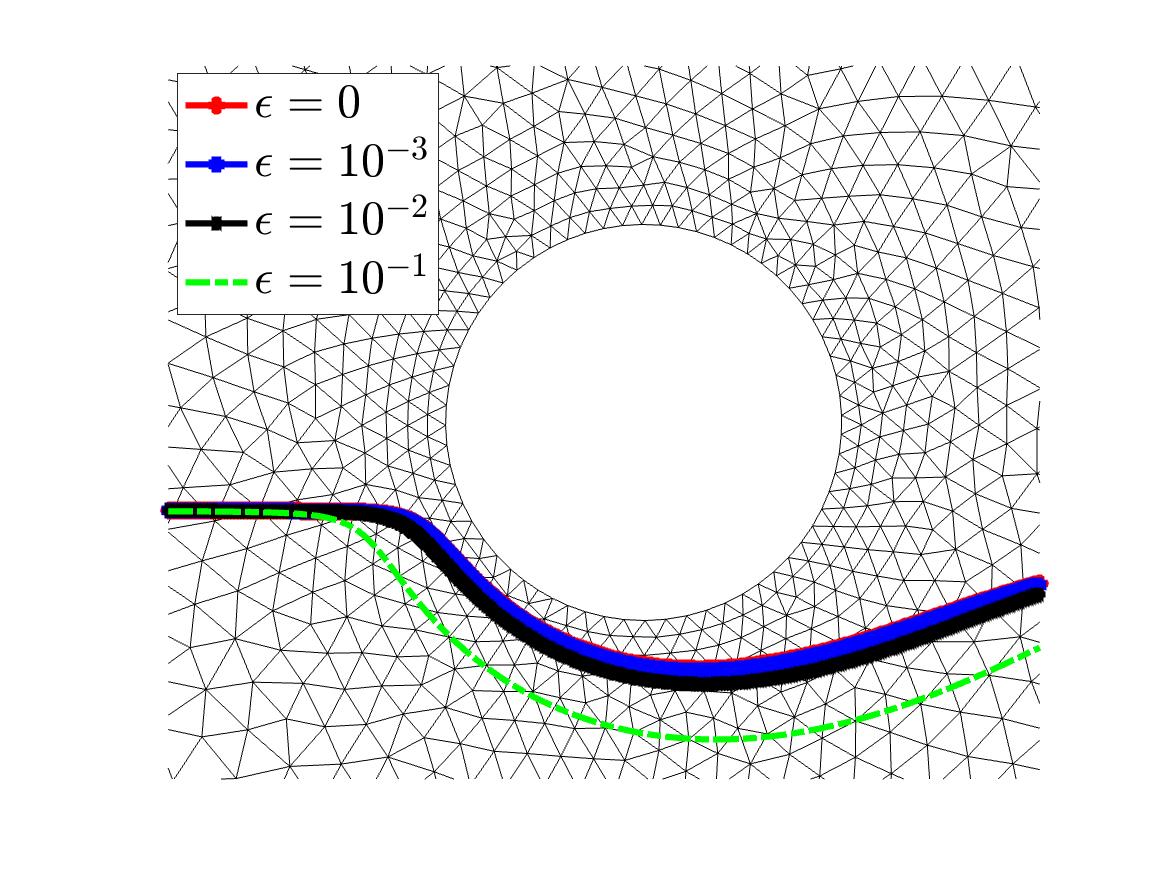}
}

\caption{transport of Gaussian distributions across a cylinder. 
(a) behavior of the $L^2$ error $\|  \rho(\cdot, t) - \rho_\infty \|_{L^1(\Omega)}$ for several choices of the regularization parameter $\epsilon$.
(b)  deformed point cloud at final time   for $\epsilon=0$ and target point cloud.
(c) behavior of select trajectories for $\epsilon=0$.
(d) behavior of one  trajectory for several values of  $\epsilon=0$.
}
\label{fig:test1_error_analysis}
\end{figure}

Figure \ref{fig:test1_error_analysis} shows the performance of our method.
Figure \ref{fig:test1_error_analysis}(a) shows the behavior of the $L^1$ error $\|\rho(\cdot, t) - \rho_\infty  \|_{L^1(\Omega)}$ for several values of $\epsilon$: as we increase $\epsilon$ in \eqref{eq:regularized_potential}, we introduce an increasing bias in the formulation that prevents convergence to the target density. Figure \ref{fig:test1_error_analysis}(b)  compares the deformed point cloud at the final time
$\{  X_j^{(K)} \}_{j}$
for $\epsilon=0$
with the ``target'' point cloud generated using $\rho_\infty$: we observe that the results are extremely satisfactory for this model problem.
 Figure \ref{fig:test1_error_analysis}(c) shows the behavior of select trajectories, while   Figure \ref{fig:test1_error_analysis}(d)   shows the behavior of a select particle   for several values of $\epsilon$: we notice that, as we increase $\epsilon$, the particle moves away from the cylinder due to the effect of the repulsive force.

\subsection{Point-set registration across a cylinder}
\label{sec:numerics_test2}

\subsubsection{Problem setting}
\begin{figure}[h!]
\centering
\subfloat[]{ 
\includegraphics[width=.33\textwidth]{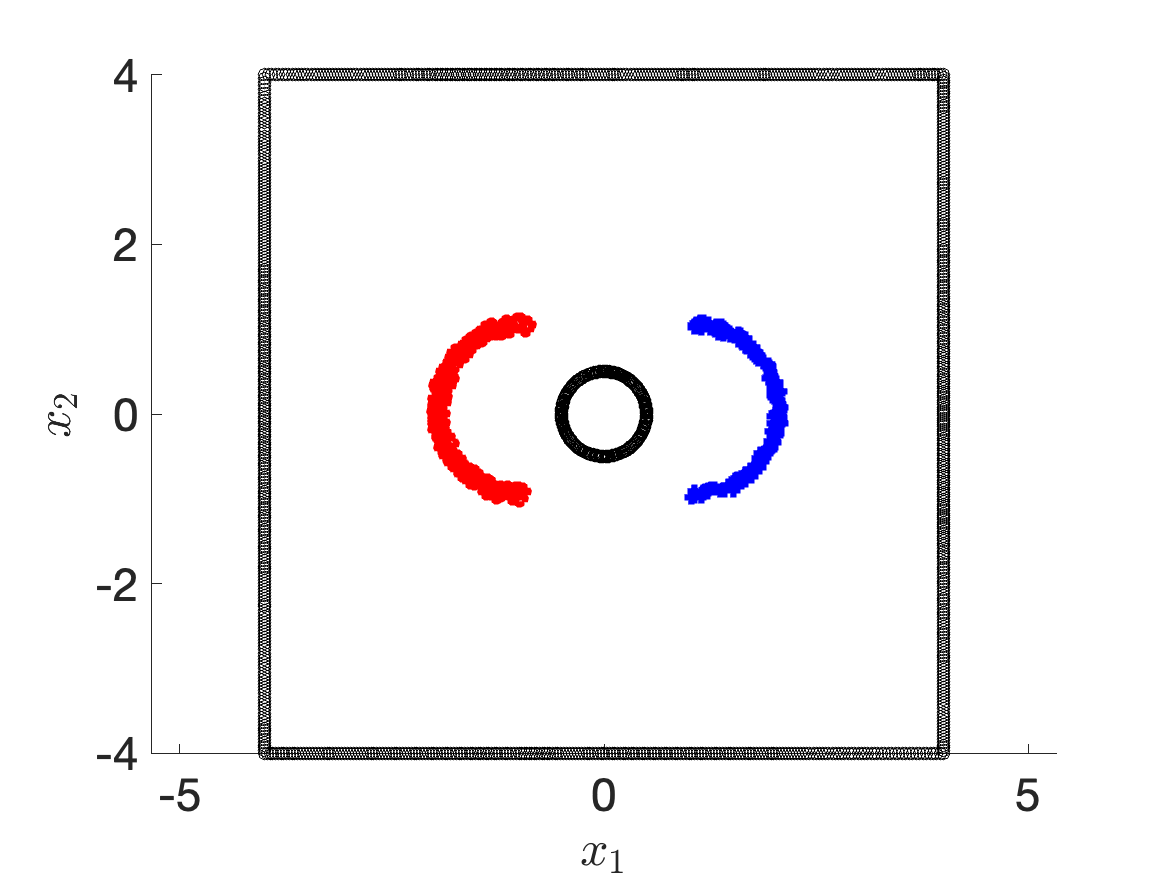}
}
~~
\subfloat[]{
\includegraphics[width=.33\textwidth]{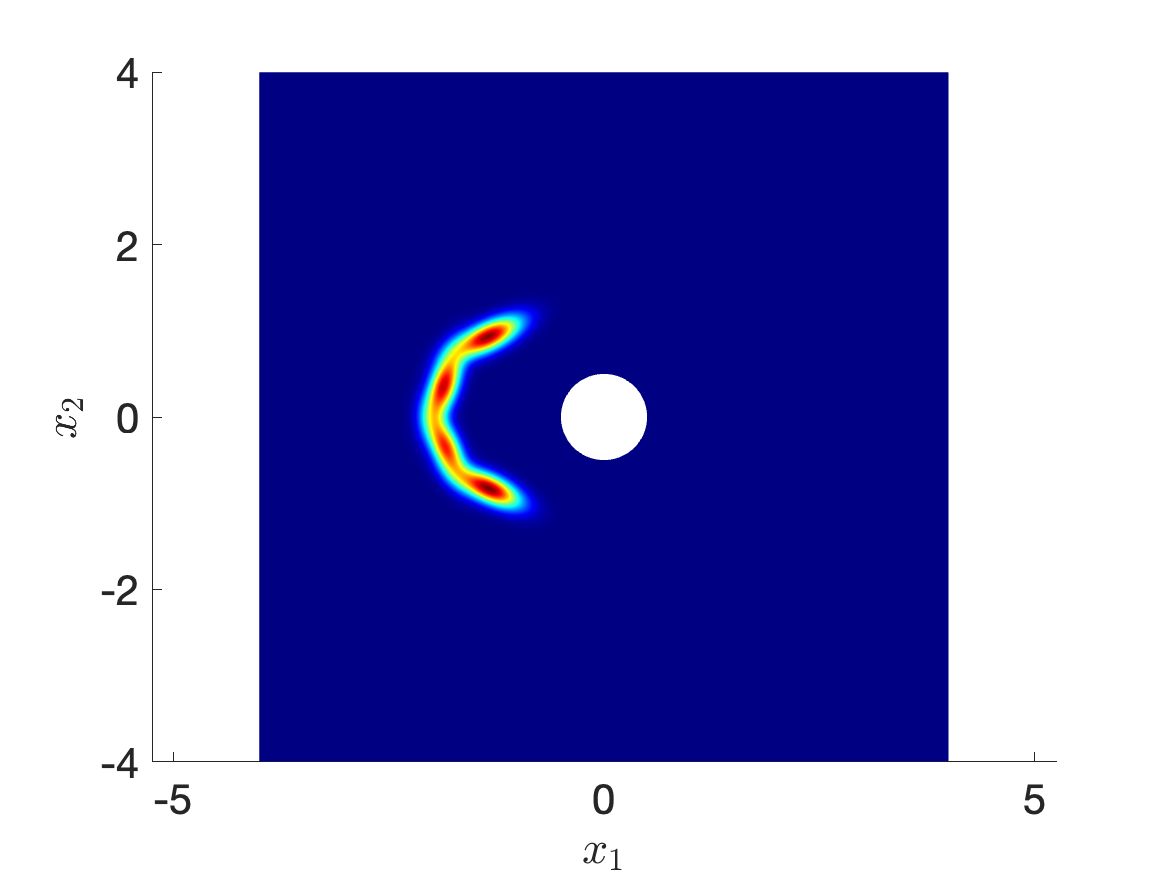}
}
~~
\subfloat[$t=0.3$]{
\includegraphics[width=0.33\textwidth]{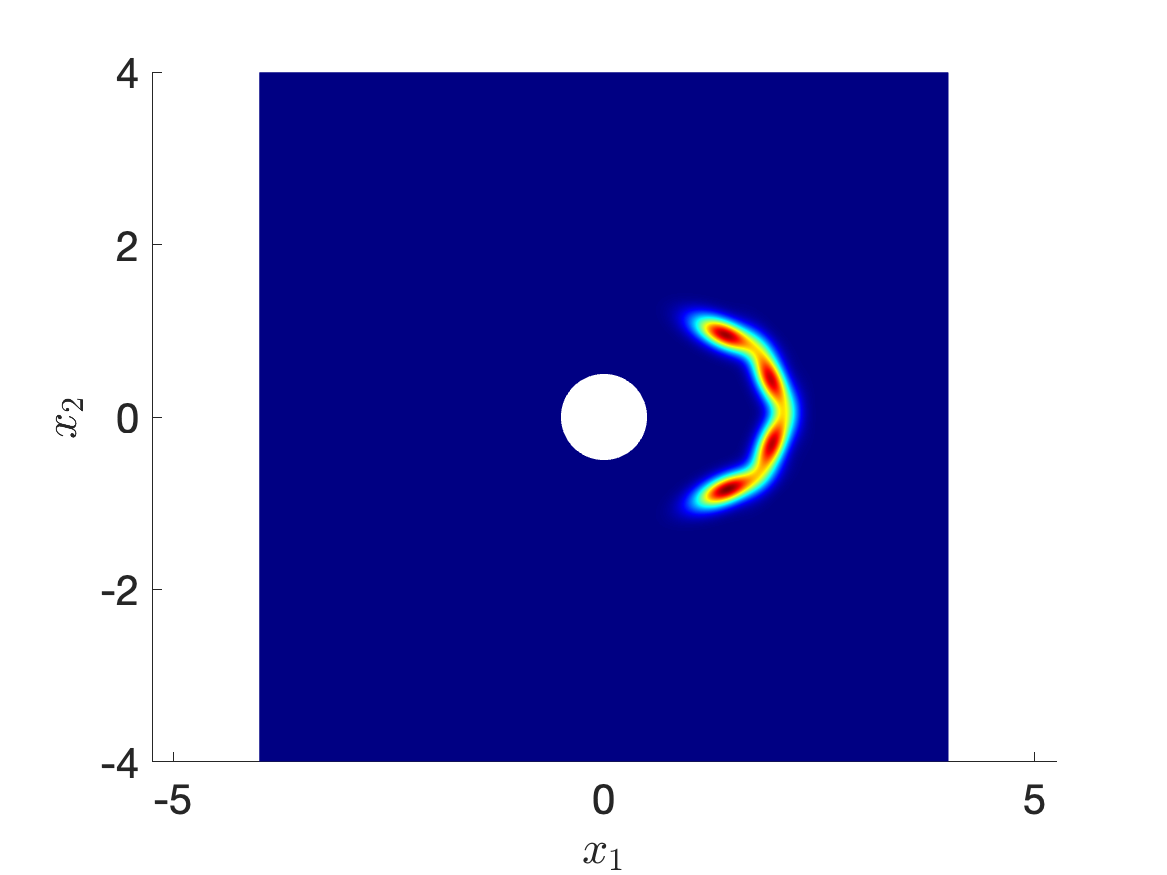}
}

\caption{point-set registration across a cylinder. 
(a) reference (red) and target point (blue) clouds.
(b) initial condition $\rho_0$.
(c) target density $\rho_\infty$.
}
\label{fig:test2_vis}
\end{figure}

We consider the problem of continuously deforming the point cloud $\{X_i^{0}\}_{i=1}^N$ to match the target point cloud $\{X_j^{\infty}\}_{j=1}^M$. We consider the same domain $\Omega$ of the previous example and we define
\begin{equation}
\label{eq:point_clouds}
\left\{
\begin{array}{ll}
\displaystyle{X_i^0 = \left[\cos(\theta_0 + \Delta \theta \frac{i-1}{N-1} ), \sin(\theta_0 + \Delta \theta \frac{i-1}{N-1} )   \right] + 0.1 \widetilde{X}_i^0},
&
i=1,\ldots,N;\\[3mm]
\displaystyle{X_j^\infty = \left[\cos(\theta_\infty + \Delta \theta \frac{j-1}{M-1} ), \sin(\theta_\infty + \Delta \theta \frac{j-1}{N-1} )   \right] + 0.1 \widetilde{X}_j^\infty}, &
j=1,\ldots,M;\\
\end{array}
\right.    
\end{equation}
where
$N=M=141$, $\theta_0=\frac{\pi}{2}$, $\theta_\infty=\frac{3\pi}{2}$,
$\Delta \theta = \pi$ and 
$\widetilde{X}_i^0,
\widetilde{X}_j^\infty 
\overset{\rm iid}{\sim} {\rm Uniform}((0,1)^2)$.
Figure \ref{fig:test2_vis}(a) shows the reference and the target point clouds.
We consider two distinct discretizations of increasing size.
The coarse discretization features a P2 FE mesh with $N_{\rm hf}=13623$ degrees of freedom and $K=5\cdot 10^3$ time steps; the fine discretization features a P3 FE mesh with $N_{\rm hf}=121425$ degrees of freedom and $K=10^4$ time steps. In both cases, we integrate the system in the time interval $(0,T=15)$ and we consider the non-uniform time grid \eqref{eq:non_uniform_timegrid}.

We rely on the unregularized potential (i.e., $\epsilon=0$ in \eqref{eq:regularized_potential}). We compare the ODE method based on forward Euler time integration with the GF method. We also run simulations for the ODE method with RK2 time integration with two time steps per time interval: since the results of the latter are comparable with the ones of the ODE method, they are not reported below. We model the point clouds \eqref{eq:point_clouds} using two GMMs with four components --- we recall  that the number of components is selected by maximizing the AIC criterion. Figures \ref{fig:test2_vis}(b) and (c) show the GMM densities that are obtained by the learning procedure.

\subsubsection{Results}

\begin{figure}[h!]
\centering
\subfloat[]{ 
\includegraphics[width=.33\textwidth]{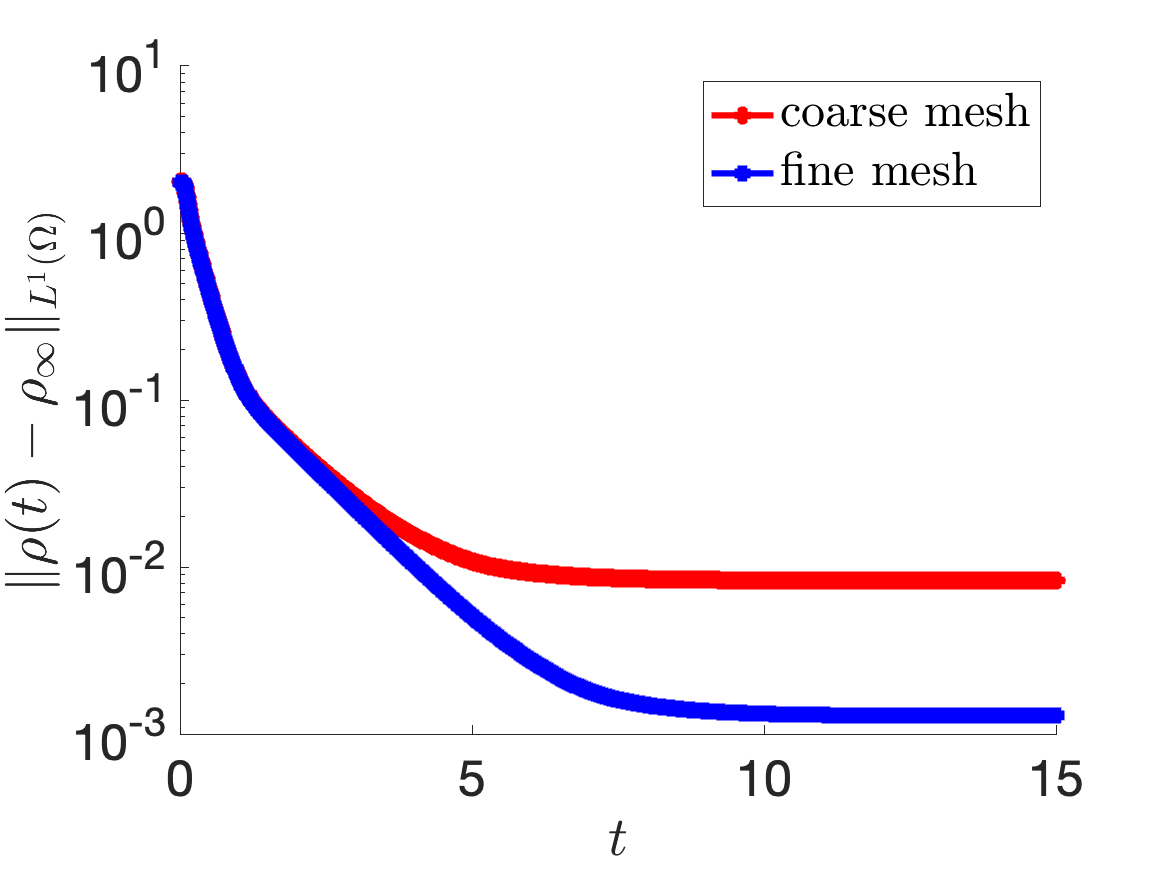}
}
~~
\subfloat[]{
\includegraphics[width=.33\textwidth]{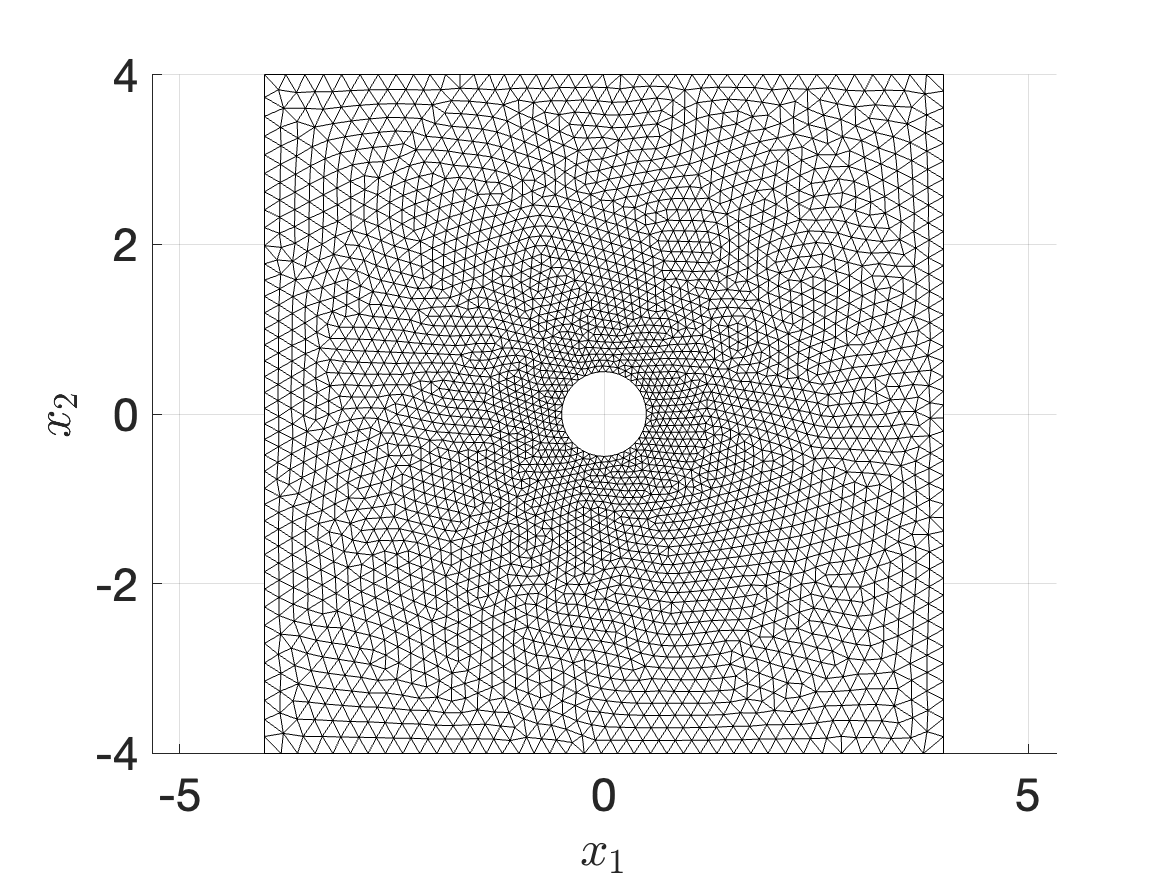}
}
~~
\subfloat[]{
\includegraphics[width=0.33\textwidth]{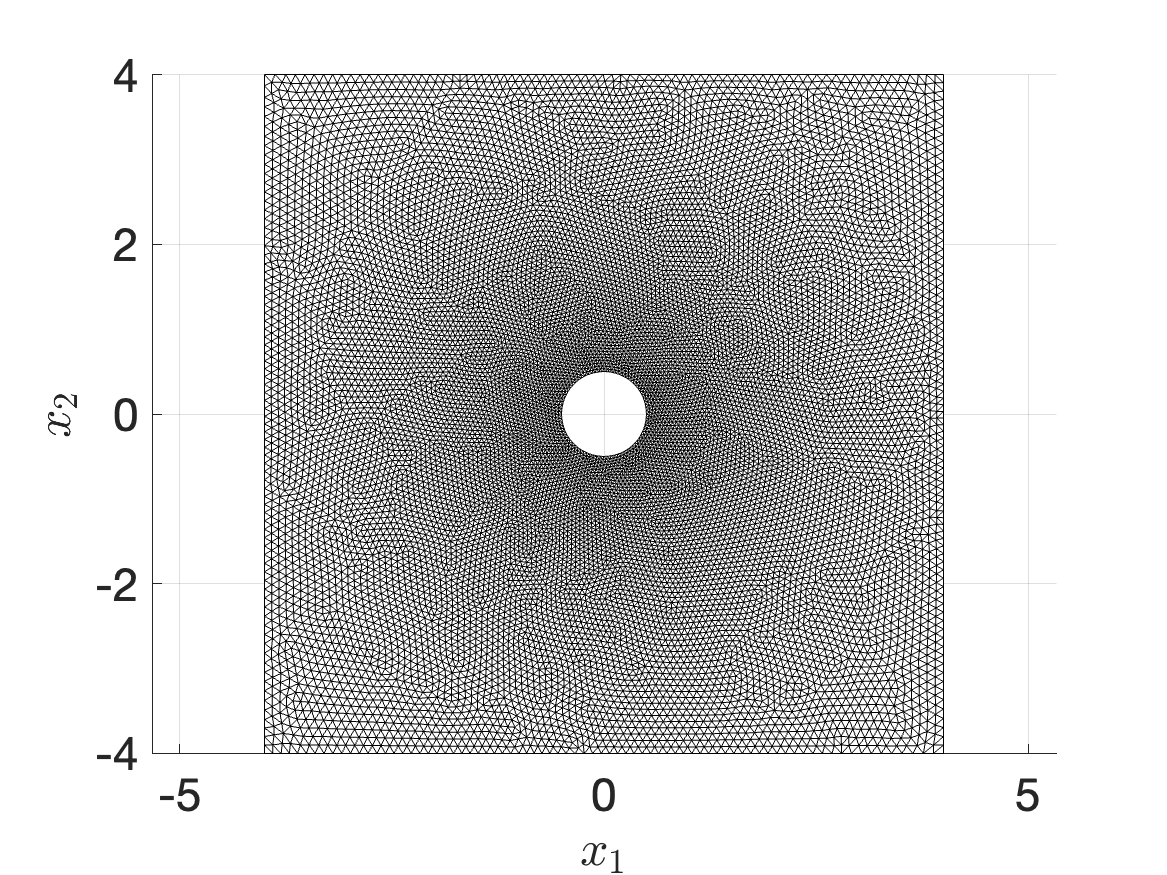}
}

\caption{point-set registration across a cylinder. 
(a) behavior of the $L^1$ error $\|  \rho(\cdot, t) - \rho_\infty \|_{L^1(\Omega)}$ for both coarse and fine meshes.
(b) coarse mesh ($N_{\rm hf}=13623$).
(c) fine mesh ($N_{\rm hf}= 121425$).
}
\label{fig:test2_error_analysis}
\end{figure}

Figure \ref{fig:test2_error_analysis}(a) shows the behavior of the 
$L^1$ error $\|  \rho(\cdot, t) - \rho_\infty \|_{L^1(\Omega)}$ over time for the two discretizations; 
Figures  \ref{fig:test2_error_analysis}(b) and (c) show the two FE meshes. As expected, the error in the target density decreases as we increase the size of the mesh.

\begin{figure}[h!]
\centering
\subfloat[$t=0.1$]{ 
\includegraphics[width=.33\textwidth]{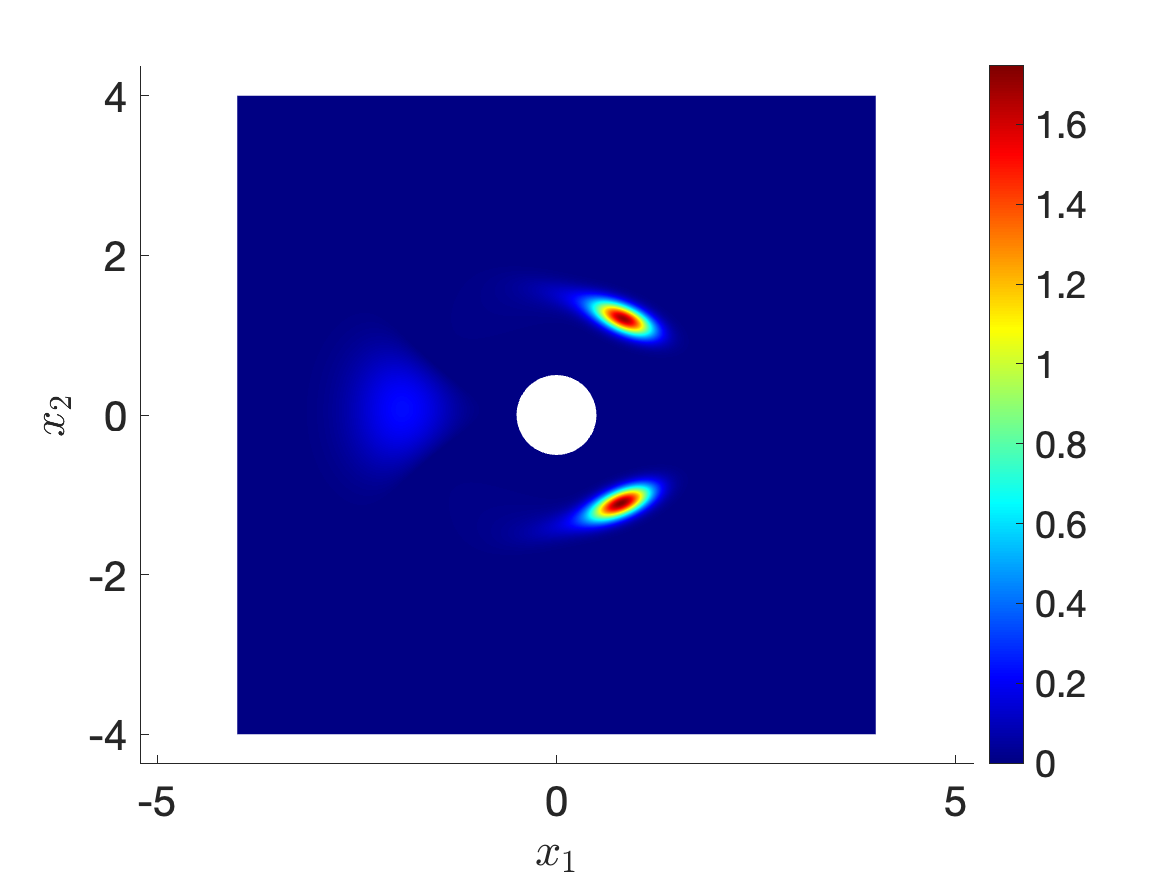}
}
~~
\subfloat[$t=0.5$]{
\includegraphics[width=.33\textwidth]{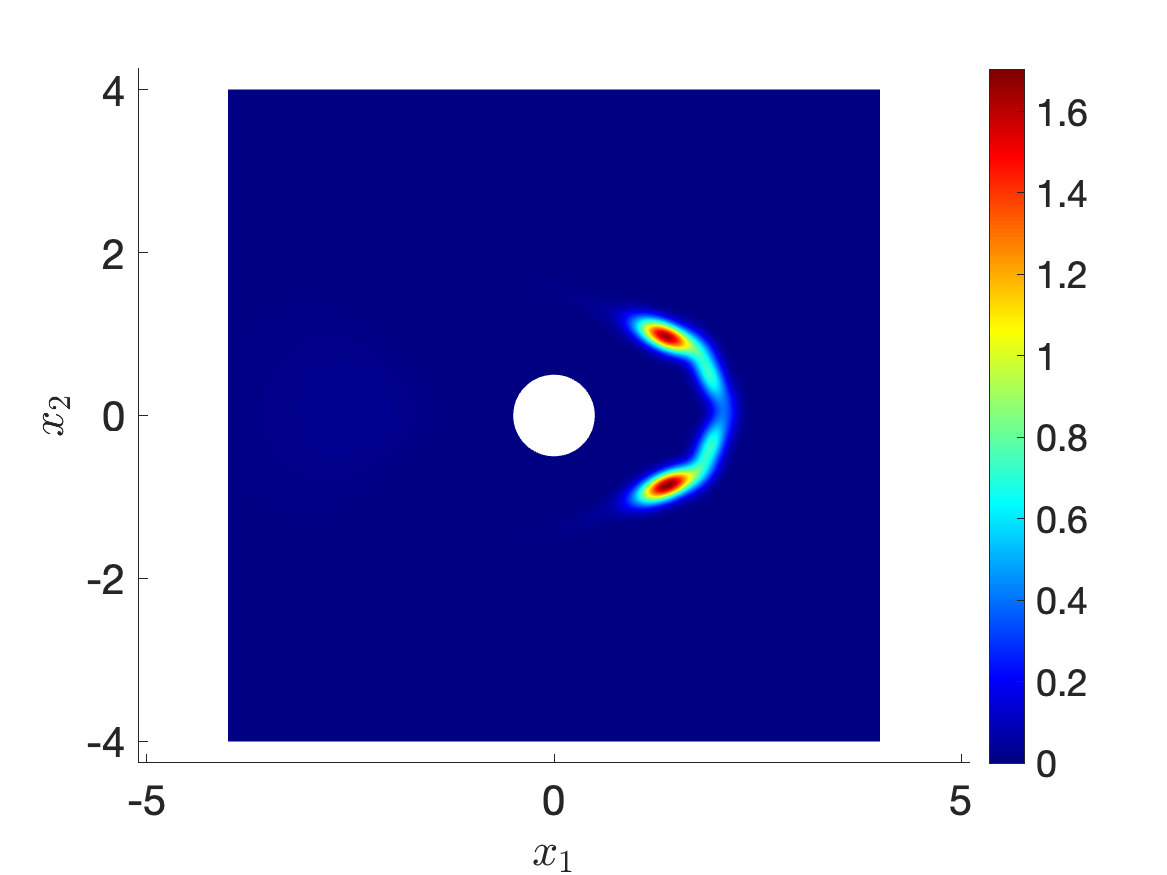}
}
~~
\subfloat[$t=1.5$]{
\includegraphics[width=0.33\textwidth]{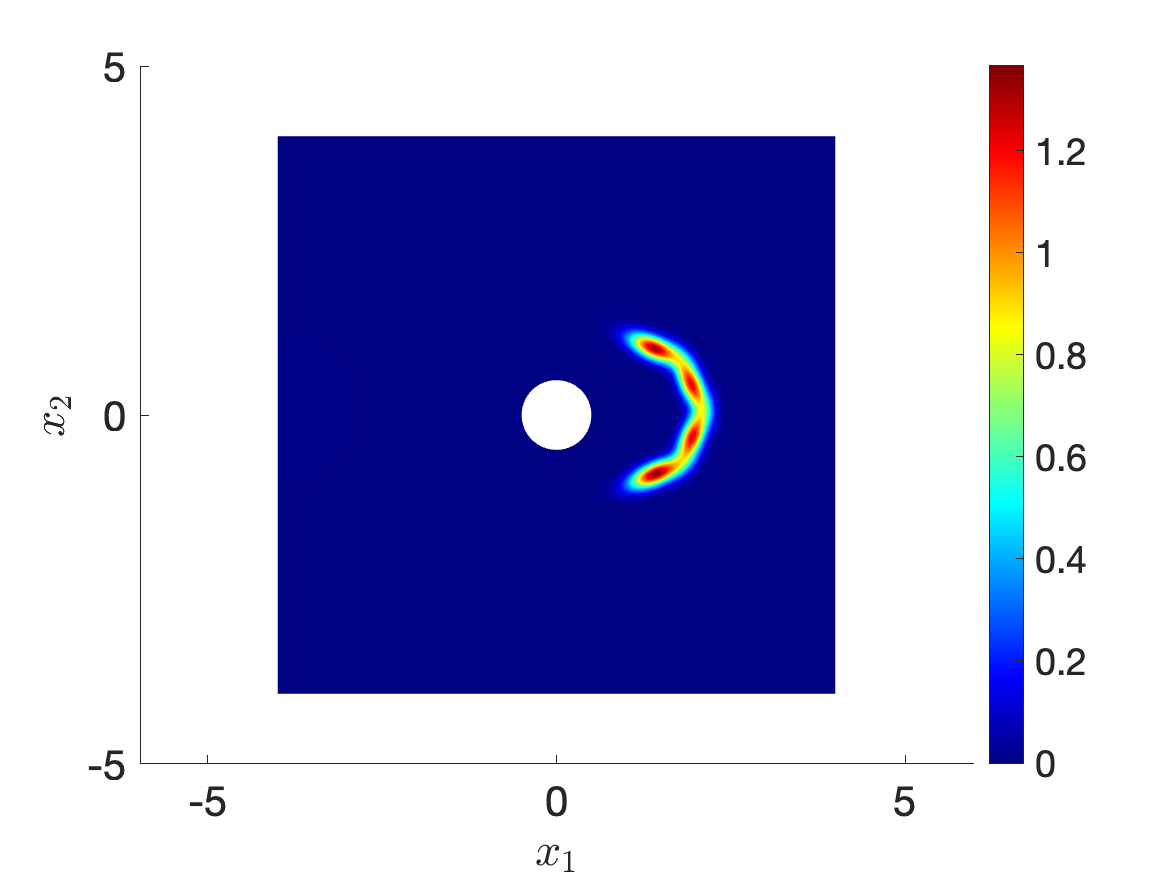}
}
\caption{point-set registration across a cylinder; density profiles for three time instants.
}
\label{fig:test2_densities}
\end{figure}

Figures \ref{fig:test2_densities}(a), (b) and (c) show the behavior of the solution to \eqref{eq:FEM_weak} for three time instants. As expected, we notice that the solution is symmetric with respect to the cylinder. We further notice that for small time steps the solution has three peaks: the majority of the mass --- corresponding to the top and the bottom peaks at $t=0.1$ --- 
rapidly goes through the cylinder, while a small portion of the mass --- corresponding to the weaker peak on the left of the cylinder at $t=0.1$--- reaches the target density later.

\begin{figure}[h!]
\centering
\subfloat[$t=0.1$]{ 
\includegraphics[width=.33\textwidth]{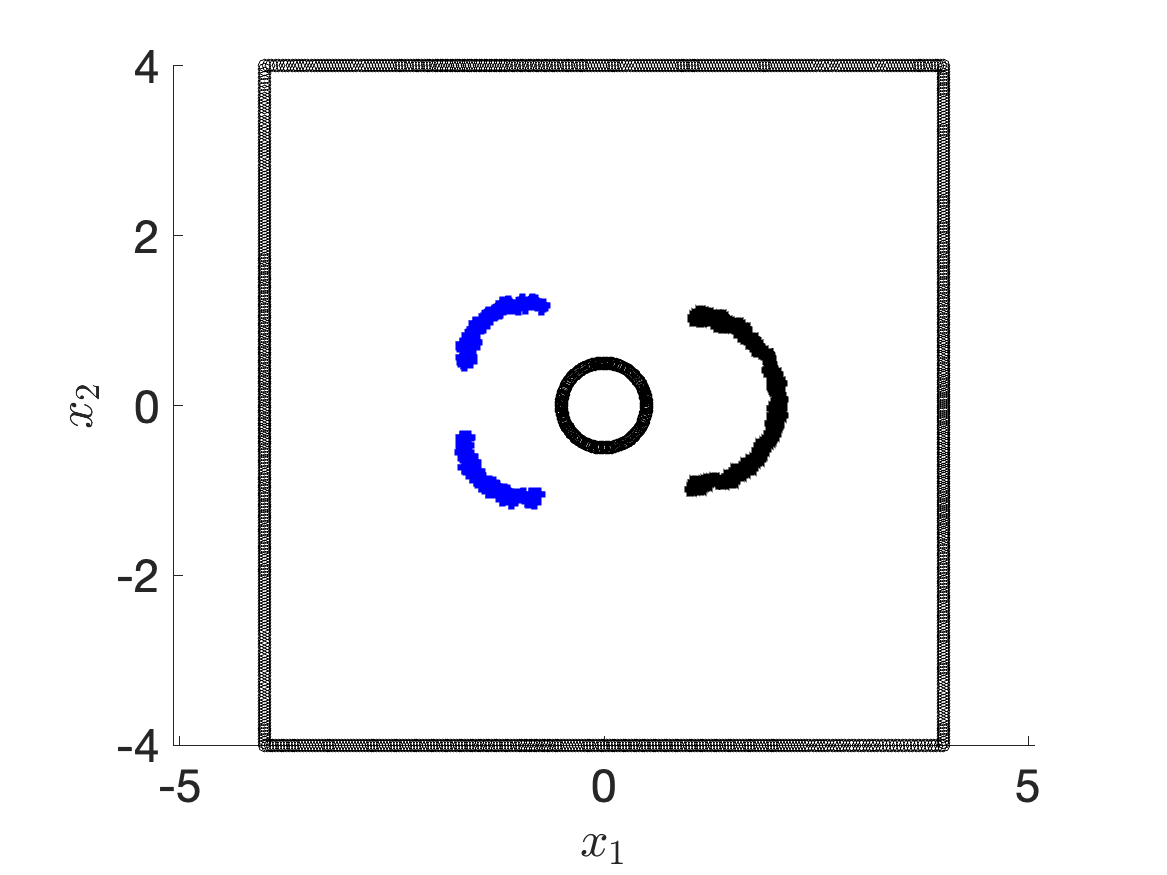}
}
~~
\subfloat[$t=0.5$]{
\includegraphics[width=.33\textwidth]{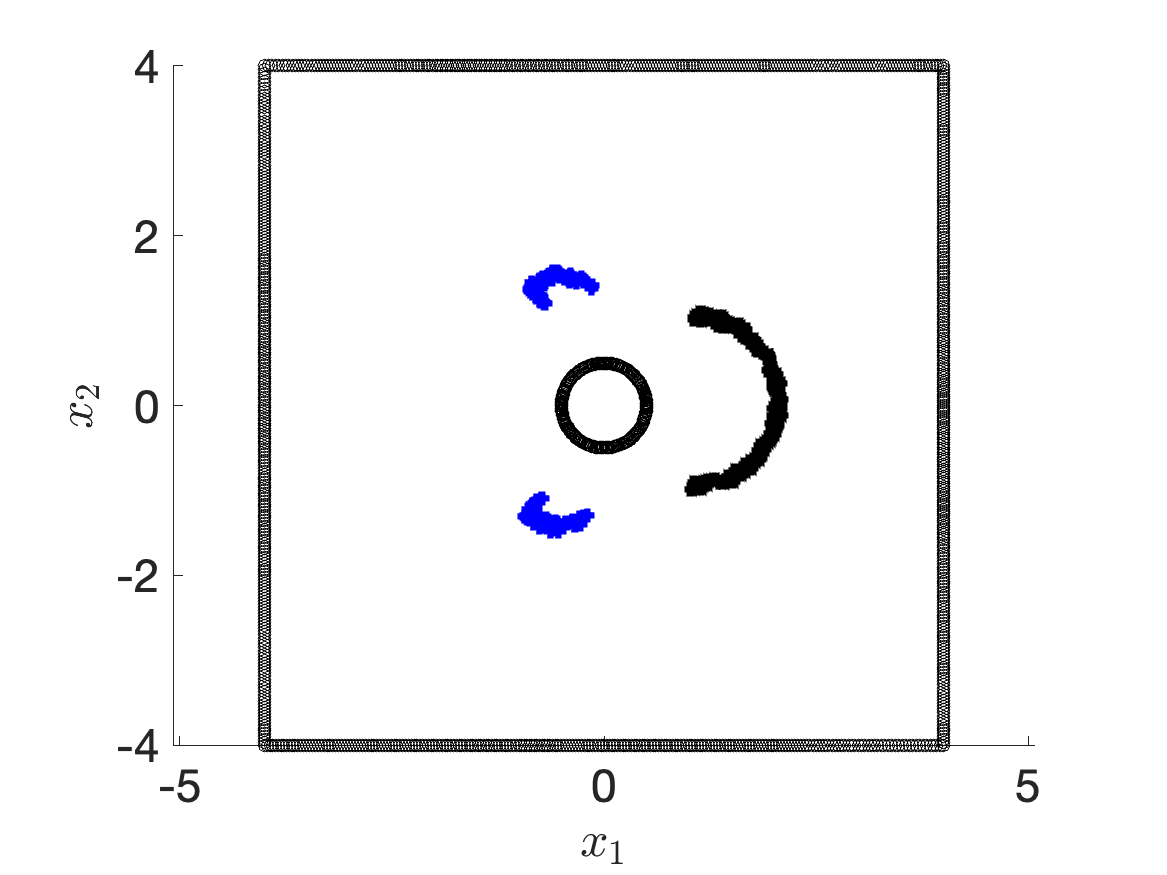}
}
~~
\subfloat[$t=1.5$]{
\includegraphics[width=0.33\textwidth]{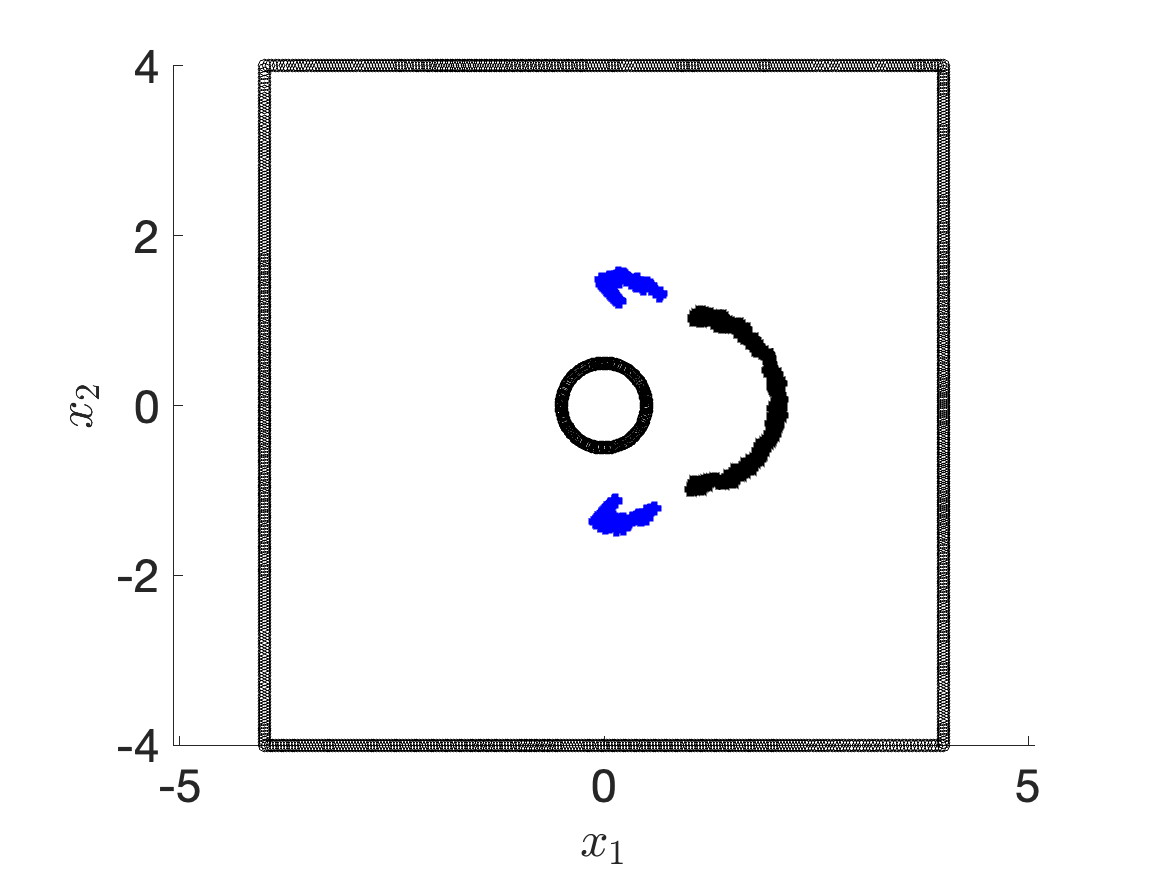}
}
\caption{point-set registration across a cylinder;  deformed point clouds $\{X_i(t)\}_{i=1}^N$ for three time instants. ODE method.
}
\label{fig:test2_point_clouds_ODE}
\end{figure}

\begin{figure}[h!]
\centering
\subfloat[$t=0.1$]{ 
\includegraphics[width=.33\textwidth]{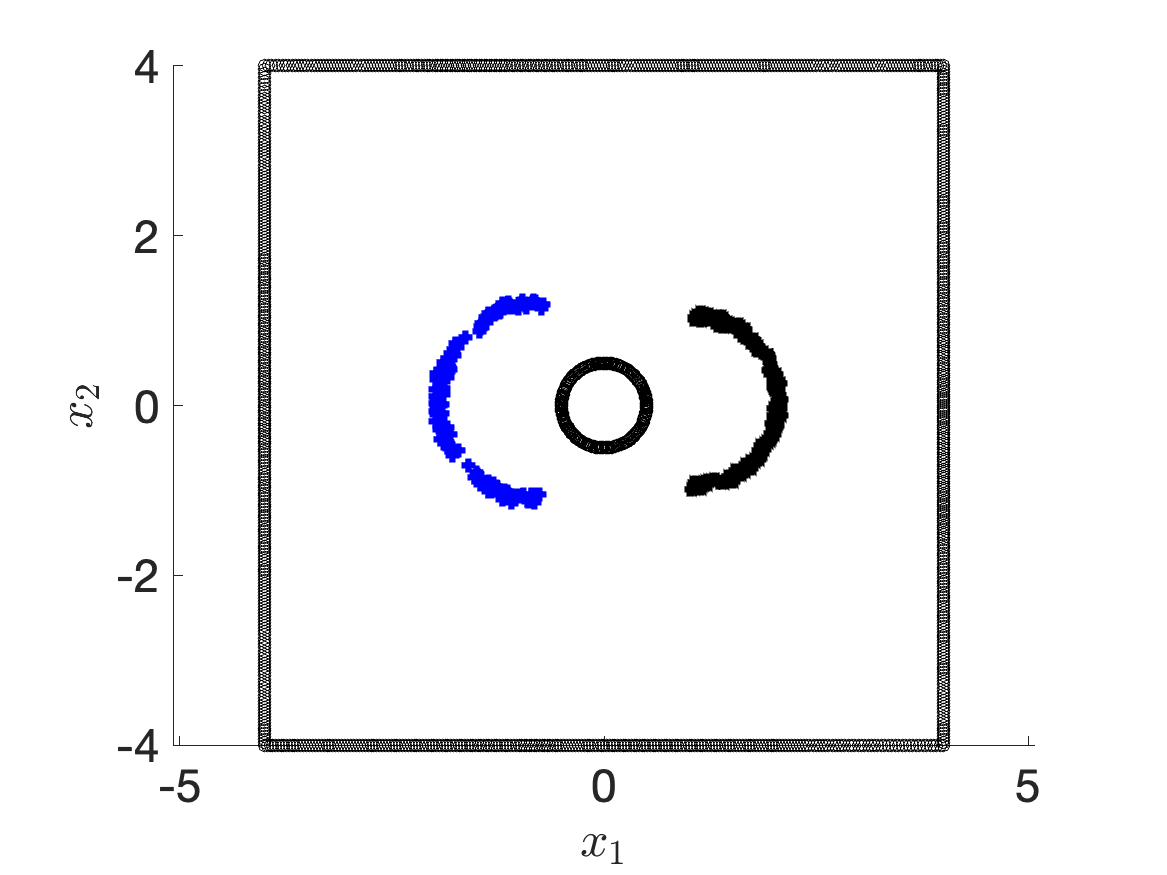}
}
~~
\subfloat[$t=0.5$]{
\includegraphics[width=.33\textwidth]{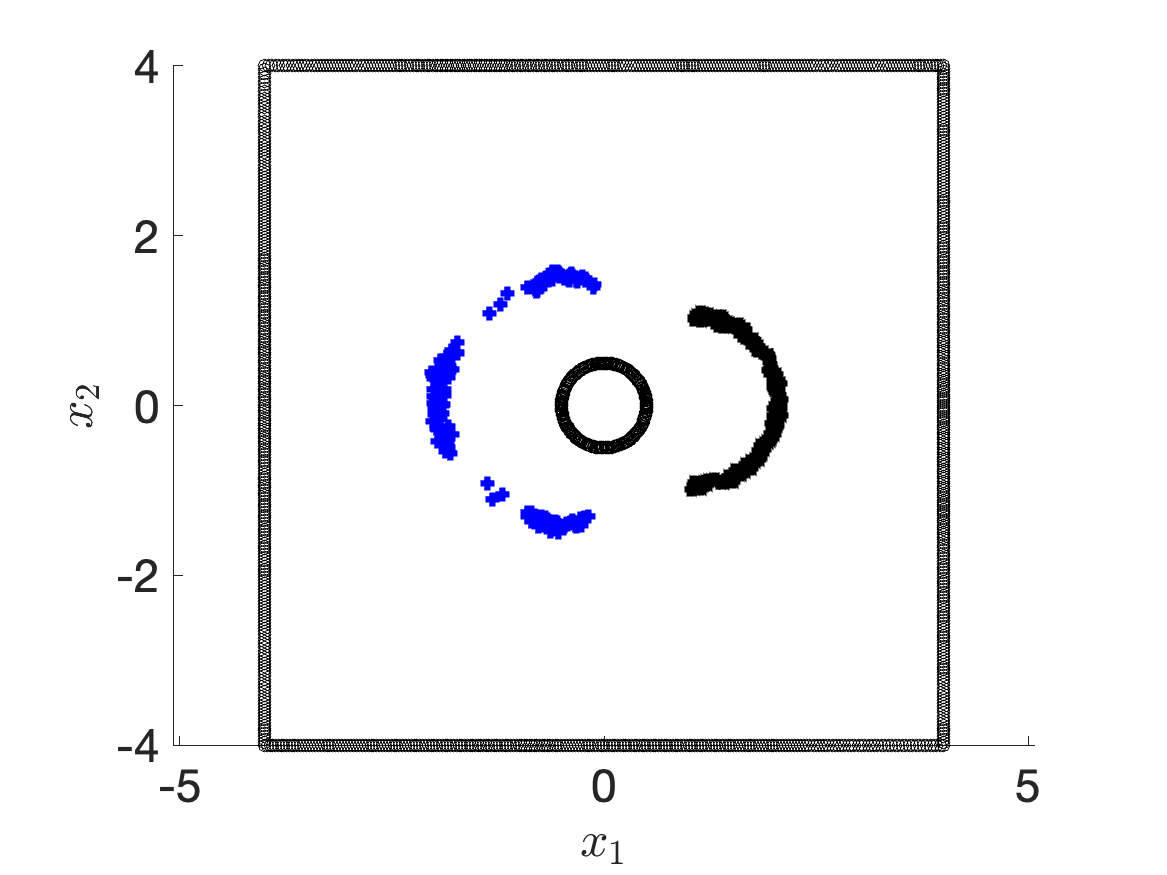}
}
~~
\subfloat[$t=1.5$]{
\includegraphics[width=0.33\textwidth]{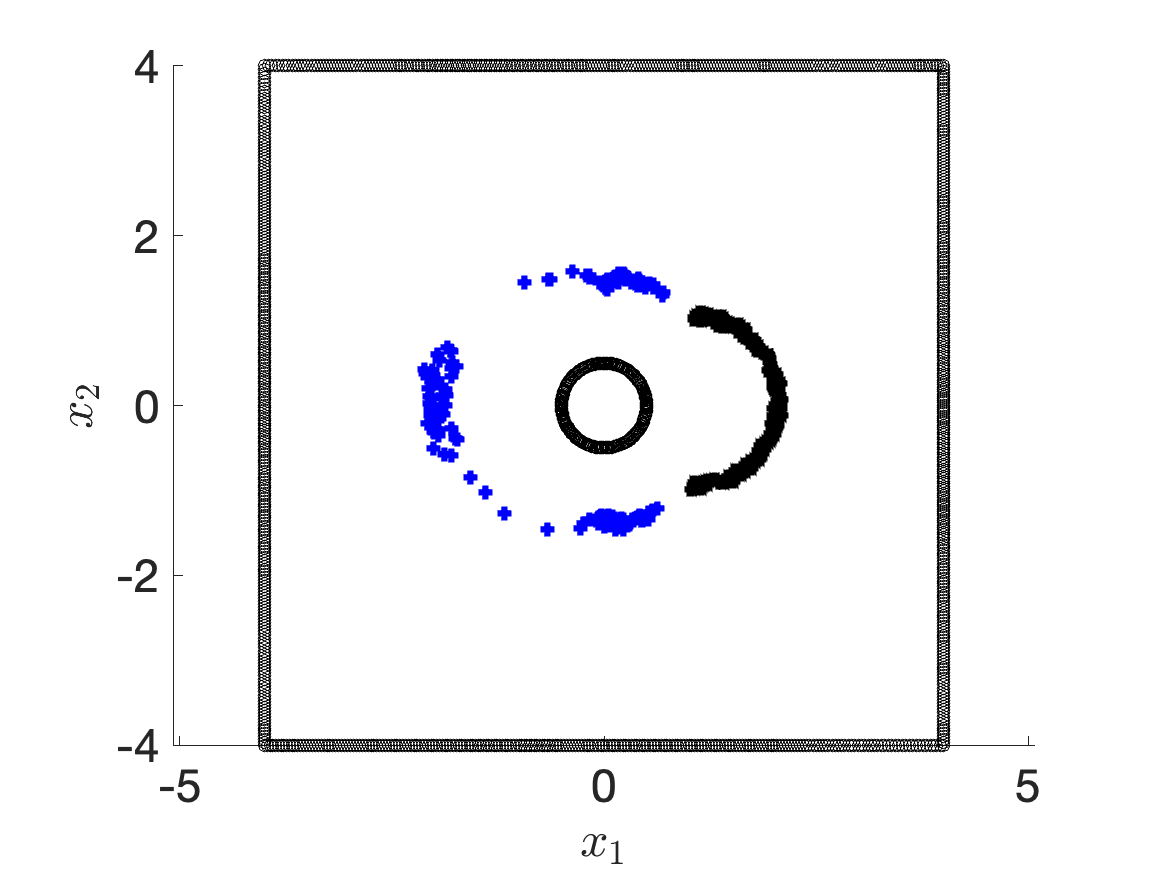}
}
\caption{point-set registration across a cylinder;  deformed point clouds $\{X_i(t)\}_{i=1}^N$ for three time instants. GF method.
}
\label{fig:test2_point_clouds_GF}
\end{figure}

Figures \ref{fig:test2_point_clouds_ODE} and \ref{fig:test2_point_clouds_GF} show the behaviors of the point clouds for the two particle integration methods discussed in section \ref{sec:methods}. 
The results correspond to the fine discretization.
Interestingly, 
the qualitative behavior of the two methods is significantly different, particularly for small values of $t$:
in the ODE method the point cloud is nearly split in two parts at $t=0.1$ and all the points proceed towards the target cloud with similar speeds; instead, in the GF method, we can clearly identify three clusters whose elements move towards the target cloud with different speeds. The latter behavior appears to be more consistent with the results for the density in Figure \ref{fig:test2_densities}. 

\begin{figure}[h!]
\centering
\subfloat[ODE]{ 
\includegraphics[width=.33\textwidth]{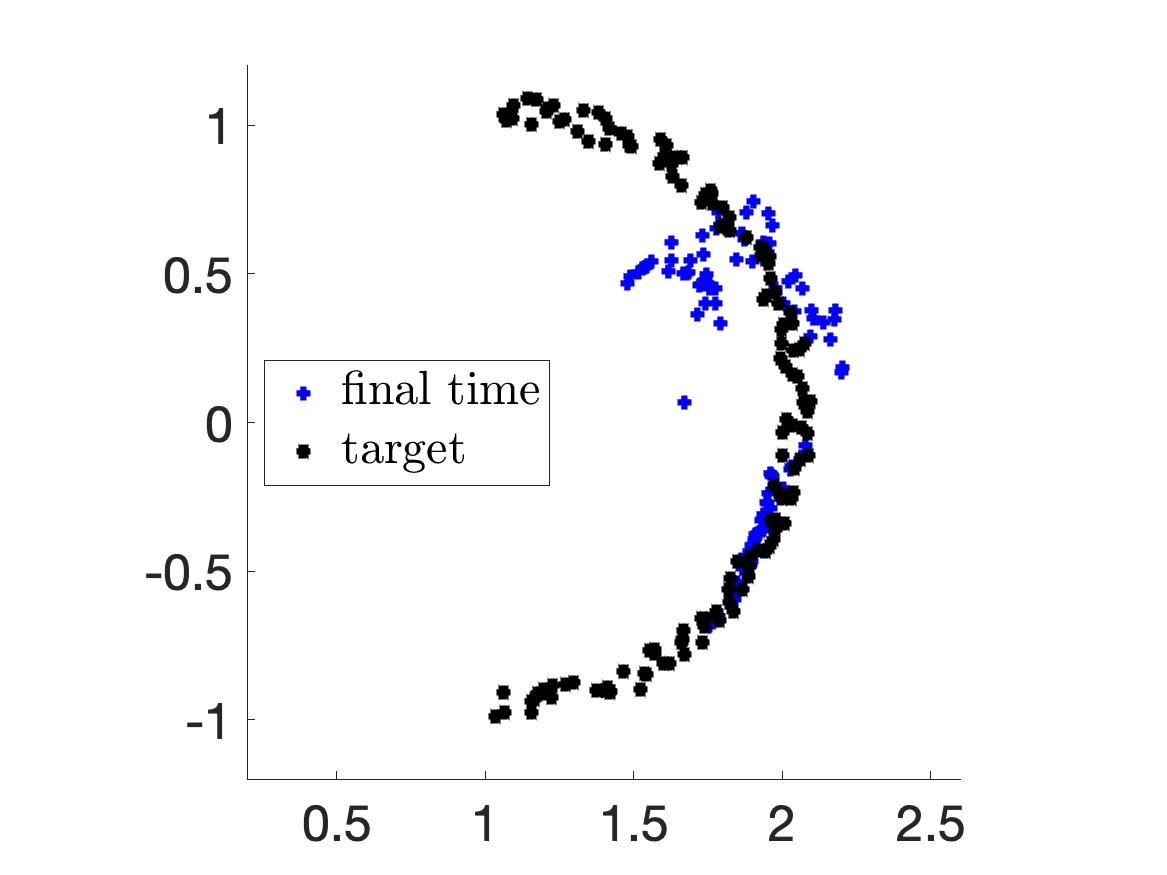}
}
~~
\subfloat[GF]{
\includegraphics[width=.33\textwidth]{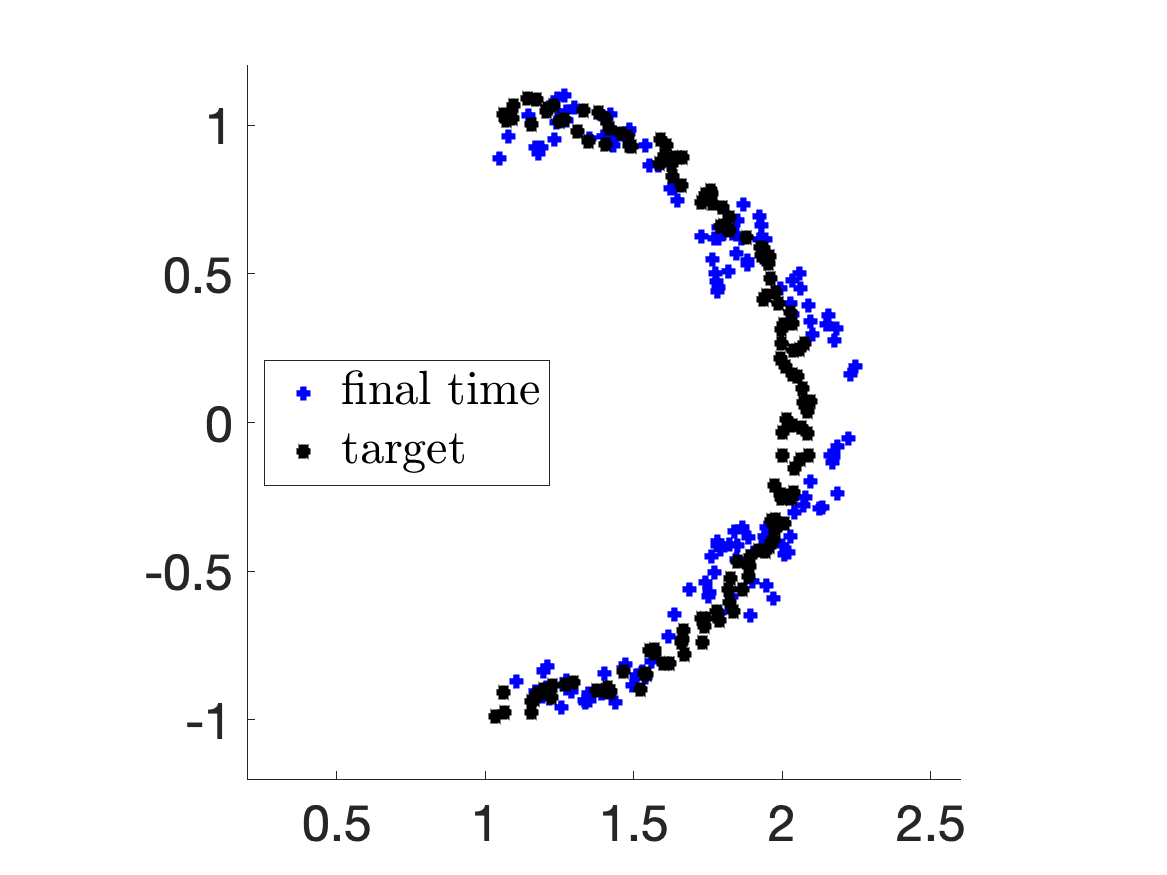}
}

\caption{point-set registration across a cylinder;  deformed point clouds $\{X_i(t)\}_{i=1}^N$  for $t=T_{\rm max}$ for the two time integration methods.}
\label{fig:test2_point_clouds_ODE_vs_GF}
\end{figure}

Figure \ref{fig:test2_point_clouds_ODE_vs_GF} compares the point clouds obtained with the two methods at the final time. We notice that the ODE method does not provide satisfactory results for this test case, while the GF method is significantly more accurate.

\section{Conclusions and perspectives}
\label{sec:conclusions}

We developed a PSR method for the interpolation of point clouds in bounded domains based on the Fokker Planck equation. Our method relies on GMMs for density estimation and on a standard (Eulerian) discretization of the FP equation for the evaluation of the density. We also proposed two distinct strategies for the estimation of the trajectories of the particles, which are both rigorously justified by the analysis of section \ref{sec:mathematical_background}. 

The numerical experiments of section \ref{sec:numerics} showed the potential of the method for PSR. The regularized potential \eqref{eq:regularized_potential} provides a simple yet effective way to control the distance of the particles from the boundaries at all times; furthermore, in both numerical experiments, the convergence to the target distribution is rapid and monotonic. 
As regards the motion of particles, the results of section \ref{sec:numerics_test2} clearly showed the superiority of the GF approach based on the solution to a  linearized OT problem over the ODE approach based on the explicit integration of the characteristics.

Despite the linearity of the equation, the solution to \eqref{eq:weak_form} using the FE method on a fixed grid is expensive and demands further investigations. We shall devise a specialized solver for the FP equation with emphasis on the accurate approximation of the particles' dynamics. 
Specifically, we plan to consider (i) adaptive time stepping schemes and (ii) adaptive spatial discretizations to reduce the overall costs: we expect to borrow ideas from deterministic particle methods for the FP equation \cite{russo1990deterministic,russo1993deterministic}.

\appendix

\section{Proofs}
\label{sec:proofs}

\subsection{Proof of Theorem \ref{th:KL}}

\begin{proof}
\emph{First statement.}
We observe that
\begin{equation}
\label{eq:silly_identity}
g(t):=t -  1 - \log (t) \geq 0 \quad
\forall \, t \in \mathbb{R}_+.
\end{equation}
To prove \eqref{eq:silly_identity}, we notice that $t=1$ is a global minimum of $g$ in $\mathbb{R}_+$ and $g(1)=0$.
If we define the set $A:=\{ x\in \Omega \,: \, \rho(x)>0 \}\subset \Omega$, we find
$$
-D_{\text{KL}}(\rho \parallel \rho_\infty) = 
\int_\Omega \rho  \, \log \left( \frac{\rho_\infty }{\rho} \right) \, dx
=
\int_A \rho  \, \log \left( \frac{\rho_\infty }{\rho} \right) \, dx;
$$
then, exploiting \eqref{eq:silly_identity}
with $t=\frac{\rho_\infty}{\rho}$, we find 
$$
-D_{\text{KL}}(\rho \parallel \rho_\infty) = 
\int_A \rho  \, \log \left( \frac{\rho_\infty }{\rho} \right)
\leq
\int_A \rho  \,  \left( \frac{\rho_\infty }{\rho}  - 1\right)
=
\int_A \rho_\infty 
-
1
\leq
\int_\Omega \rho_\infty 
-
1
=0,
$$
which proves $D_{\text{KL}}(\rho \parallel \rho_\infty)\geq 0$.

Let $\rho \in \mathcal{P}_2(\Omega)$ satisfy 
$D_{\text{KL}}(\rho \parallel \rho_\infty)=0$. If we denote by $\phi: t \mapsto - \log (t)$ and $f: x \mapsto  \frac{\rho(x)}{\rho_\infty(x)}$, we find
$$
\phi \left( \int_{\Omega} f(x) \rho(x) \, dx \right)
= \phi \left( 1 \right) = 0,
\quad
\int_\Omega \phi \circ f (x) \rho(x) \, dx = 
- D_{\text{KL}}(\rho \parallel \rho_\infty) = 0.
$$
Notice that 
$\phi \left( \int_{\Omega} f(x) \rho(x) \, dx \right) = 
\int_\Omega \phi \circ f (x) \rho(x) \, dx$:
since $f$ is convex, recalling Jensen's inequality, we must have that $f$ is constant, which implies that 
$\rho=  \rho_\infty$ a.s..
\end{proof}

\begin{proof}
\emph{Second statement}.
It suffices to show an example for which 
$D_{\text{KL}}(\rho \parallel \rho_\infty) \neq
D_{\text{KL}}(\rho_\infty \parallel \rho)$.
Towards this end,  consider $\Omega=(-1,1)$ and 
$$
\rho(x) \equiv \frac{1}{2},
\quad
\rho_\infty(x) =
\left\{
\begin{array}{ll}
0.1     &  x\in (-1,0)\\
0.9     &  x \in (0,1)
\end{array}
\right.
$$
By straightforward calculations, we find
$D_{\text{KL}}(\rho \parallel \rho_\infty) =
0.5108$ and 
$D_{\text{KL}}(\rho_\infty \parallel \rho) =
0.3681$.
\end{proof}

\begin{proof}
\emph{Third statement}.
 We define $\rho_\lambda := \lambda \rho_0 + (1-\lambda) \rho_1$. Since the function
$t\mapsto t \log(t/c)$ is strictly convex for any $c>0$, we have
$$
\rho_\lambda(x) \log \left(  \frac{\rho_\lambda(x)}{\rho_\infty(x)} \right)
\leq
\lambda \rho_0(x) \log \left(  \frac{\rho_0(x)}{\rho_\infty(x)} \right)
+
(1- \lambda) \rho_1(x) \log \left(  \frac{\rho_1(x)}{\rho_\infty(x)} \right),
\quad
\forall \, x\in \Omega; \quad
\lambda\in (0,1);
$$
furthermore, the latter holds with equality if and only if $\rho_0(x)=\rho_1(x)$.
Recalling the expression of the KL divergence, we obtain the desired result.
\end{proof}

\begin{proof}
\emph{Fourth statement}.
See \cite[Theorem 3]{gilardoni2010pinsker}.
\end{proof}

\subsection{Proof of Theorem \ref{th:FPequation}}

\begin{enumerate}
\item 
The existence and uniqueness of the solution to \eqref{eq:weak_form} can be studied using the theory of abstract parabolic problems. We refer to
\cite[Chapter 10]{salsa2015partial} (cf. Theorem 10.11) for a detailed proof.
\item 
Exploiting 
\cite[Theorem 10.11]{salsa2015partial}, we find that 
$\rho(\cdot, t) \in L^2(\Omega)$ for all $t>0$; it hence  suffices to show that $\int_\Omega \rho(x,t) \, dx = 1$ for all $t>0$ and that $\rho(x,t) \geq 0$ 
for a.e. $x\in \Omega$ and $t>0$.
If we consider  $v=1$ in \eqref{eq:weak_form}, we find
$$
0=
\int_\Omega  \partial_t \rho(x,t) \, dx
=
\frac{d}{dt}
\int_\Omega  \rho(x,t) \, dx  \quad \Rightarrow \quad
\int_\Omega  \rho(x,t) \, dx = 
\int_\Omega  \rho_0(x) \, dx = 1;
$$
which proves conservation of mass. The proof of positivity is significantly more technical and is provided in 
\cite[Corollary 4.3]{ouhabaz2009analysis}.
\item 
Recalling notation introduced in \eqref{eq:KLdivergence_revisited},
since $\rho \geq 0$ a.e., 
we find that 
$$
\frac{d}{dt}D_{\text{KL}}(\rho(\cdot,t) \parallel \rho_\infty) 
=
- \int_\Omega \rho(x,t) \big| \nabla F'[\rho(x,t)] \big|^2 \, dx;
$$
Provided that 
$c:=\inf_{x\in \Omega} \rho_\infty(x) $ is strictly positive  and
$\rho_0 \in L^2(\Omega)$, we have that 
$D_{\text{KL}}(\rho_0 \parallel \rho_\infty)<\infty$:
we have indeed 
$$
D_{\text{KL}}(\rho_0 \parallel \rho_\infty)
\leq
\int_\Omega \rho_0  \log  \left( \frac{\rho_0 }{c} \right)
=
\int_\Omega \rho_0  \log \rho_0  \, dx
-
\log(c) 
\leq
\int_\Omega (\rho_0 ^2 + \rho_0  ) \, dx
-
\log(c) 
=
\| \rho_0   \|_{L^2(\Omega)}^2
+1 - \log(c),
$$
where in the second inequality we used the identity $t\log (t) \leq t^2+t$ that is valid for all $t\in \mathbb{R}_+$.
Therefore, since 
$\rho(\cdot, t)\in \mathcal{P}_2(\Omega)$ 
for all $t>0$,
 we find that the KL divergence is monotonically-decreasing and uniformly bounded. Exploiting the Pinsker's inequality, we also find that $t\mapsto \|  \rho(\cdot, t) \|_{L^1(\Omega)}$ is uniformly bounded.
Then, we can apply Theorem 2.3 in \cite{alves2024strong} to conclude that $\rho(\cdot,t) \to \rho_\infty$ in $L^1$ for $t\to \infty$. 
\item 
The result exploits the Lagrangian interpretation of the FP equation, which is extensively discussed in 
\cite[Chapter 4]{santambrogio2015optimal} --- the precise result is provided in
Theorem 4.4 of \cite{santambrogio2015optimal}.
\end{enumerate}

\bibliographystyle{abbrv}
\bibliography{all_refs}

\begin{thebibliography}{10}

\bibitem{alves2024strong}
N.~J. Alves, J.~Skrzeczkowski, and A.~E. Tzavaras.
\newblock Strong convergence of sequences with vanishing relative entropy.
\newblock {\em arXiv preprint arXiv:2409.16892}, 2024.

\bibitem{brenier1991polar}
Y.~Brenier.
\newblock Polar factorization and monotone rearrangement of vector-valued
  functions.
\newblock {\em Communications on pure and applied mathematics}, 44(4):375--417,
  1991.

\bibitem{brooks1982streamline}
A.~N. Brooks and T.~J. Hughes.
\newblock Streamline upwind/{P}etrov-{G}alerkin formulations for convection
  dominated flows with particular emphasis on the incompressible
  {N}avier-{S}tokes equations.
\newblock {\em Computer methods in applied mechanics and engineering},
  32(1-3):199--259, 1982.

\bibitem{cucchiara2024model}
S.~Cucchiara, A.~Iollo, T.~Taddei, and H.~Telib.
\newblock Model order reduction by convex displacement interpolation.
\newblock {\em Journal of Computational Physics}, 514:113230, 2024.

\bibitem{fox2012tutorial}
C.~W. Fox and S.~J. Roberts.
\newblock A tutorial on variational {B}ayesian inference.
\newblock {\em Artificial intelligence review}, 38:85--95, 2012.

\bibitem{gilardoni2010pinsker}
G.~L. Gilardoni.
\newblock On {P}insker's and {V}ajda's type inequalities for {C}sisz{\'a}r's $
  f $-divergences.
\newblock {\em IEEE Transactions on Information Theory}, 56(11):5377--5386,
  2010.

\bibitem{gray2011entropy}
R.~M. Gray.
\newblock {\em Entropy and information theory}.
\newblock Springer Science \& Business Media, 2011.

\bibitem{hastie2009elements}
T.~Hastie, R.~Tibshirani, and J.~Friedman.
\newblock {\em The Elements of Statistical Learning}.
\newblock Springer Series in Statistics. Springer New York Inc., New York, NY,
  USA, 2009.

\bibitem{iollo2014advection}
A.~Iollo and D.~Lombardi.
\newblock Advection modes by optimal mass transfer.
\newblock {\em Physical Review E}, 89(2):022923, 2014.

\bibitem{iollo2022mapping}
A.~Iollo and T.~Taddei.
\newblock Mapping of coherent structures in parameterized flows by learning
  optimal transportation with {G}aussian models.
\newblock {\em Journal of Computational Physics}, 471:111671, 2022.

\bibitem{jordan1998variational}
R.~Jordan, D.~Kinderlehrer, and F.~Otto.
\newblock The variational formulation of the {F}okker--{P}lanck equation.
\newblock {\em SIAM journal on mathematical analysis}, 29(1):1--17, 1998.

\bibitem{kullback1997information}
S.~Kullback.
\newblock {\em Information theory and statistics}.
\newblock Courier Corporation, 1997.

\bibitem{mozolevski2007hp}
I.~Mozolevski, E.~S{\"u}li, and P.~R. B{\"o}sing.
\newblock hp-version a priori error analysis of interior penalty discontinuous
  {G}alerkin finite element approximations to the biharmonic equation.
\newblock {\em Journal of Scientific Computing}, 30(3):465--491, 2007.

\bibitem{ouhabaz2009analysis}
E.-M. Ouhabaz.
\newblock {\em Analysis of heat equations on domains}.
\newblock Princeton University Press, 2009.

\bibitem{peyre2019computational}
G.~Peyr{\'e}, M.~Cuturi, et~al.
\newblock Computational optimal transport: With applications to data science.
\newblock {\em Foundations and Trends{\textregistered} in Machine Learning},
  11(5-6):355--607, 2019.

\bibitem{russo1990deterministic}
G.~Russo.
\newblock Deterministic diffusion of particles.
\newblock {\em Communications on Pure and Applied Mathematics}, 43(6):697--733,
  1990.

\bibitem{russo1993deterministic}
G.~Russo.
\newblock A deterministic vortex method for the {N}avier-{S}tokes equations.
\newblock {\em Journal of Computational Physics}, 108(1):84--94, 1993.

\bibitem{salsa2015partial}
S.~Salsa.
\newblock {\em Partial differential equations in action}, volume~1.
\newblock Springer, third edition, 2016.

\bibitem{santambrogio2015optimal}
F.~Santambrogio.
\newblock Optimal transport for applied mathematicians.
\newblock {\em Birk{\"a}user, NY}, 55(58-63):94, 2015.

\bibitem{song2019generative}
Y.~Song and S.~Ermon.
\newblock Generative modeling by estimating gradients of the data distribution.
\newblock {\em Advances in neural information processing systems}, 32, 2019.

\bibitem{taddei2020registration}
T.~Taddei.
\newblock A registration method for model order reduction: data compression and
  geometry reduction.
\newblock {\em SIAM Journal on Scientific Computing}, 42(2):A997--A1027, 2020.

\bibitem{van1992stochastic}
N.~G. Van~Kampen.
\newblock {\em Stochastic processes in physics and chemistry}, volume~1.
\newblock Elsevier, 1992.

\end{thebibliography}
  
\end{document}